\documentclass[oneside,english]{amsart}
\usepackage[T1]{fontenc}
\usepackage[latin1]{inputenc}
\usepackage{amssymb}

\makeatletter
 \theoremstyle{plain}    
 \newtheorem{thm}{Theorem}[section]
 \numberwithin{equation}{section} %% Comment out for sequentially-numbered
 \numberwithin{figure}{section} %% Comment out for sequentially-numbered
 \theoremstyle{plain}
 \theoremstyle{plain}    
 \newtheorem{lem}[thm]{Lemma} %%Delete [thm] to re-start numbering
 \theoremstyle{plain}    
 \newtheorem{prop}[thm]{Proposition} %%Delete [thm] to re-start numbering
 \theoremstyle{definition}
 \newtheorem{defn}[thm]{Definition}
 \theoremstyle{remark}
 \newtheorem{rem}[thm]{Remark}
 \theoremstyle{plain}    
 \newtheorem{cor}[thm]{Corollary} %%Delete [thm] to re-start numbering
 \theoremstyle{definition}
  \newtheorem{example}[thm]{Example}

\usepackage[all]{xy}

\makeatletter
\newcommand{\xyR}[1]{%
\makeatletter
\xydef@\xymatrixrowsep@{#1}
\makeatother
} % end of \xyR

\makeatletter
\newcommand{\xyC}[1]{%
\makeatletter
\xydef@\xymatrixcolsep@{#1}
\makeatother
}% end of \xyC

\usepackage{babel}
\makeatother
\begin{document}

\title{On inner Kan complexes in the category of dendroidal sets}

\author{I. Moerdijk, I. Weiss}

\begin{abstract}
The category of dendroidal sets is an extension of that of simplicial
sets, suitable for defining nerves of operads rather than just of
categories. In this paper, we prove some basic properties of inner
Kan complexes in the category of dendroidal sets. In particular, we
extend fundamental results of Boardman and Vogt, of Cordier and Porter,
and of Joyal to dendroidal sets. 
\end{abstract}
\maketitle

\section{Introduction}

This paper is a companion to our paper \cite{DenSet}, where we introduced
the category of dendroidal sets, and explained some of its applications.
The main goal of the present paper is to introduce a notion of ''inner
Kan complex'' for dendroidal sets, and to prove some of the fundamental
properties of these inner Kan complexes. 

Dendroidal sets provide a generalization of simplicial sets. There
is an embedding $i_{!}:sSet\rightarrow dSet$ along which many properties
of simplicial sets can be extended to dendroidal sets. For example,
the category of dendroidal sets carries a closed symmetric monoidal
structure, which extends the Cartesian structure of simplicial sets
(but is not itself Cartesian). This monoidal structure is closely
related to the Boardman-Vogt tensor product of operads (\cite{BoardmanVogt,Fiedor,Dunn}).
It plays a central role in various uses of dendroidal sets, for example
in the construction of homotopy coherent nerves of operads, and in
the definition of weak higher categories given in \cite{DenSet}.

We recall that a simplicial set is said to satisfy the Kan condition,
or, to be a Kan complex, if every horn $\Lambda^{i}[n]\rightarrow X$
for $0\le i\le n$ has a filler. Boardman and Vogt \cite{BoardmanVogt}
study this filler condition for $0<i<n$, and refer to it as the restricted
Kan condition. Recently, Joyal \cite{JoyalBook} has been studying simplicial sets satisfying this condition, under the name of quasi-categories.  The horns $\Lambda^{i}[n]$ for $i$ different from
$0$ and $n$ are called inner horns, and we shall call a simplicial
set satisfying this restricted Kan condition an inner Kan complex.
Note that the nerve of any category is an inner Kan complex.

In this paper, we will define inner horns and inner Kan complexes
for dendroidal sets, in such a way that a simplicial set $X$ is an
inner Kan complex iff $i_{!}(X)$ is a dendroidal inner Kan complex.
The dendroidal nerve of any operad provides an example of such a dendroidal
inner Kan complex. We will prove several fundamental properties of
dendroidal inner Kan complexes. Our main result is that the closed
monoidal structure on dendroidal sets has the property that for any
two dendroidal sets $X$ and $Y$, the internal Hom $\underline{Hom}(X,Y)$
is an inner Kan complex whenever  $X$ is normal and $Y$ is inner Kan (Theorem 9.1 below). We will
also show that the homotopy coherent nerve of a topological operad
is an inner Kan complex. (A more general statement for operads in
monoidal model categories is given in Theorem 7.1 below.) These results
specialise to known results for simplicial sets, and provide new proofs
of these. Indeed, Cordier and Porter \cite{CordierPorter} prove that
the homotopy coherent nerve of a locally fibrant simplicial category
is an inner Kan complex. And Joyal \cite{JoyalBook} proves for any
two simplicial sets $X$ and $Y$ that $\underline{Hom}(X,Y)$ is
an inner Kan complex (quasi-category) whenever $Y$ is (see also \cite{Joshua}).
The latter result plays a fundamental role in Joyal's proof of the
existence of a model structure on simplicial sets in which the inner
Kan complexes are the fibrant objects. We expect our result for dendroidal
sets to play a similar role in establishing an analogous model structure
on the category of dendroidal sets. 

In \cite{DenSet}, we introduced a Grothendieck construction (homotopy
colimit) for diagrams of dendroidal sets, necessary for our definition
of weak higher categories. In Section 8 of this paper, we will prove
that this dendroidal Grothendieck construction yields an inner Kan
complex when applied to a diagram of inner Kan complexes. In addition,
in Section 6, we will give an explicit description of the operad generated
by an inner Kan complex, modelled on the one given by Boardman and
Vogt \cite{BoardmanVogt} in the simplicial case. We then use this
description to prove that a dendroidal set satisfies the unique filler
condition for inner horns iff it is the nerve of an operad.

\section{The category of dendroidal sets}

The notion of dendroidal set was introduced in \cite{DenSet}. We
briefly recall the relevant definitions here. 

To begin with, we introduce a category $\Omega$ whose objects are
finite rooted trees. If we think of a tree as a graph, and call a
vertex unary if it has only one edge attached to it, then our trees
$T$ are equipped with a distinguished unary vertex $o$ called the
output, and a set of unary vertices $I$ (not containing $o$) called
the set of inputs. When drawing such a tree, it is common to orient
the tree ''towards the output'' drawn at the bottom, and delete
the designated output and input vertices from the picture. Thus, in
the tree $T$, 

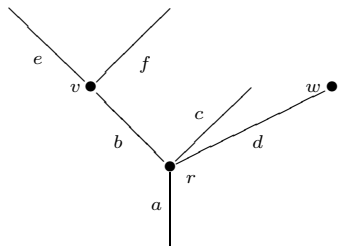
\[
\xymatrix{*{\,}\ar@{-}[dr]_{e} &  & *{\,}\ar@{-}[dl]^{f}\\
\,\ar@{}[r]|{\,\,\,\,\,\,\,\,\,\,\,\,\,\, v} & *{\bullet}\ar@{-}[dr]_{b} &  & *{\,}\ar@{-}[dl]_{c}\ar@{}[r]|{\,\,\,\,\,\,\,\,\,\,\,\, w} & *{\bullet}\ar@{-}[dll]^{d}\\
 &  & *{\bullet}\ar@{-}[d]_{a} & \,\ar@{}[l]^{r\,\,\,\,\,\,\,\,\,\,\,}\\
 &  & *{\,}}
\]
the output vertex at the edge $a$ has been deleted, as have the input
vertices at $e,f,c$. This tree $T$ now has three (remaining) vertices,
three \emph{input edges} attached to the three deleted input vertices,
and one \emph{output edge} (attached to the deleted output vertex).
These input and output edges are called \emph{outer edges} (the output
edge is also called the \emph{root}, while the input edges are also
called \emph{leaves}), the others ($b$ and $d$ in the picture) are
called \emph{inner} edges. From now on, we will not mention the input
and output vertices anymore, and ''vertex'' will always refer to
a remaining vertex. 

Attached to each such vertex in the tree, there will be one designated
edge pointing towards the root; the other edges attached to this vertex
are called the input edges of that vertex, and their number is called
the \emph{valence} of the vertex. So in the tree $T$ pictured above,
the vertex $r$ has valence three and the vertex $w$ has valence
zero. The tree with just one edge is now drawn as \[
\xymatrix{*{}\ar@{-}[d]\\
*{}}
\]
and referred to as $\eta$, or sometimes as $\eta_{e}$ if we want
to name its edge $e$. The linear tree, with one input edge and one
output edge and $n$ vertices, is denoted $i[n]$. It has $n+1$ edges
which we usually number from input to output as $0,1,\cdots,n$. Here
is a picture of $i[2]$:
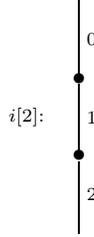
\[
\xymatrix{ & *{}\ar@{-}[d]^{0}\\
\ar@{}[d]^{i[2]:} & *{\bullet}\ar@{-}[d]^{1}\\
*{} & *{\bullet}\ar@{-}[d]^{2}\\
 & *{}}
\]
Each tree $T$ defines a coloured operad (see \cite{ColWCons} for
a definition of coloured operads) which we denoted $\Omega(T)$ in
\cite{DenSet}. The colours of this operad are the edges of the tree,
and the operations are \emph{generated} by the vertices of the tree.
A planar representation of the tree gives a specific set of generators.
For example, for the tree $T$ pictured above, $\Omega(T)$ has six
colours, $a,b,\cdots,f$. A choice of generating operations is $r\in\Omega(T)(b,c,d;a)$,
$w\in\Omega(T)(-;d)$ and $v\in\Omega(T)(e,f;b)$. The other operations
are units such as $1_{b}\in\Omega(T)(b;b)$, compositions such as
$r\circ_{1}v\in\Omega(T)(e,f,c,d;a)$, and permutations such as $r\cdot\tau\in\Omega(T)(c,b,d;a)$.
Note that the \emph{same} tree $T$ can be given a different planar
structure, e.g. 
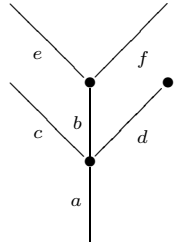
\[
\xymatrix{ & *{}\ar@{-}[dr]_{e} &  & *{}\ar@{-}[dl]^{f}\\
 & *{}\ar@{-}[dr]_{c} & *{\bullet}\ar@{-}[d]_{b} & *{\bullet}\ar@{-}[dl]^{d}\\
 &  & *{\bullet}\ar@{-}[d]_{a}\\
 &  & *{}}
\]
 which defines the \emph{same} operad $\Omega(T)$ but suggests a
different choice of generators ($r\cdot\tau$ rather than $r$). 

The category $\Omega$ is now defined as the category having these
trees $T$ (with designated output and inputs) as objects, and as
arrows $T\rightarrow T'$ the maps of coloured operads $\Omega(T)\rightarrow\Omega(T')$.
(Note that every such map sends colours to colours, i.e., edges of
$T$ to edges of $T'$, and is in fact completely determined by this). 

The category of \emph{dendroidal sets} is the category of functors
$X:\Omega^{op}\rightarrow Set$ and natural transformations between
them. We will denote this category by $dSet$. 

For a dendroidal set $X$ and a tree $T$, we will usually write $X_{T}$
for $X(T)$, and call an element of the set $X_{T}$ a \emph{dendrex}
of $X$ of ''shape'' $T$. 

The linear trees $i[n]$ for $n\ge0$ define a functor (a full embedding)\[
i:\Delta\rightarrow\Omega\]
 from the standard simplicial category, and hence by composition a
functor \[
i^{*}:dSet\rightarrow sSet,\]
 from dendroidal sets to the category $sSet$ of simplicial sets.

By Kan extension, this functor has both a left and a right adjoint,
denoted $i_{!}$ and $i_{*}:sSet\rightarrow dSet$, respectively.
The functor $i_{!}$ is ''extension by zero''; for a simplicial
set $X$,\[
i_{!}(X)_{T}=\left\{ \begin{array}{cc}
X_{n}, & \textrm{if $T\cong i[n]$ for some $n\ge0$}\\
\phi, & \textrm{otherwise}\end{array}\right.\]
 This defines a full embedding $i_{!}:sSet\rightarrow dSet$, from
simplicial sets into dendroidal sets. 

Each coloured operad $\mathcal{P}$ defines a dendroidal set $N_{d}(\mathcal{P})$,
its \emph{dendroidal nerve}, by\[
N_{d}(\mathcal{P})_{T}=Hom(\Omega(T),\mathcal{P}),\]
\emph{Hom} denoting the set of arrows in the category of operads.
If $\mathcal{P}$ is itself an operad of the form $\Omega(S)$, given
by an object $S$ of $\Omega$, then $N_{d}(\mathcal{P})$ is the
representable dendroidal set given by $S$, which we will denote by
$\Omega[S]$; i.e.,\[
N_{d}\Omega(S)=\Omega[S]\]
by definition. The functor $N_{d}$ from operads to dendroidal sets
is fully faithful, and has a left adjoint which we will denote by
\[
\tau_{d}:dSet\rightarrow Operad.\]
For a dendroidal set $X$, we refer to $\tau_{d}(X)$ as \emph{the
operad generated} by $X$.

We also recall from \cite{DenSet} that the Cartesian structure on
$sSet$ extends to a (non-Cartesian) closed symmetric monoidal structure
$\otimes$ on $dSet$. This structure is completely determined by
the identity\[
\Omega[S]\otimes\Omega[T]=N_{d}(\Omega(S)\otimes_{BV}\Omega(T))\]
where $\otimes_{BV}$ denotes the Boardman-Vogt tensor product of
(coloured) operads; see \cite{BoardmanVogt}. The corresponding internal
Hom is then determined by the Yoneda lemma, as\[
\underline{Hom}(X,Y)_{T}=Hom_{dSet}(\Omega[T]\otimes X,Y).\]
We will come back to the monoidal structure in more detail in Section
9.

\section{Faces and degeneracies}

Exactly as for $\Delta$, the maps in $\Omega$ are generated by special
kinds of maps. 

(i) Given a tree $T$ and a vertex $v\in T$ of valence 1, there is
a tree $T'$, obtained from $T$ by deleting the vertex $v$ and merging
the two edges $e_{1}$ and $e_{2}$ on either side of $v$ into one
new edge $e$. There is an operad map, i.e. an arrow $\sigma_{v}:T\rightarrow T'$
in $\Omega$, which sends $v$ to the unit $1_{e}$. For example:\[
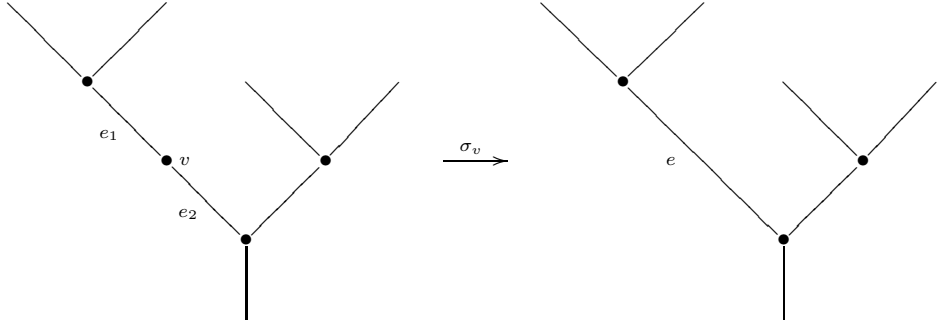
\begin{array}{ccc}
\xymatrix{*{\,}\ar@{-}[dr] &  & *{\,}\ar@{-}[dl]\\
 & *{\bullet}\ar@{-}[dr]_{e_{1}} &  & *{\,}\ar@{-}[dr] &  & *{\,}\ar@{-}[dl]\\
 &  & *{\bullet}\ar@{-}[dr]_{e_{2}}\ar@{}|{\,\,\,\,\,\,\,\,\,\, v} &  & *{\bullet}\ar@{-}[dl]\\
 &  &  & *{\bullet}\ar@{-}[d]\\
 &  &  & *{\,}}
 & \xymatrix{\\\\\ar[r]^{\sigma_{v}} & *{}}
 & \xymatrix{*{\,}\ar@{-}[dr] &  & *{\,}\ar@{-}[dl]\\
 & *{\bullet}\ar@{-}[ddrr] &  & *{\,}\ar@{-}[dr] &  & *{\,}\ar@{-}[dl]\\
 & \,\ar@{}[r]|{\,\,\,\, e} &  &  & *{\bullet}\ar@{-}[dl]\\
 &  &  & *{\bullet}\ar@{-}[d]\\
 &  &  & *{\,}}
\end{array}\]

An arrow in $\Omega$ of this kind will be called a \emph{degeneracy
(map).}

(ii) Given a tree $T$, and a vertex $v$ in $T$ with exactly one
inner edge attached to it, one can obtain a new tree $T/v$ by deleting
$v$ and all the external edges attached to it. The operad $\Omega(T/v)$
associated to $T/v$ is simply a suboperad of the one associated to
$T$, and this inclusion of operads defines an arrow in $\Omega$
denoted \[
\partial_{v}:T/v\rightarrow T.\]
An arrow in $\Omega$ of this kind is called an \emph{outer face (map).}
For example\emph{\[
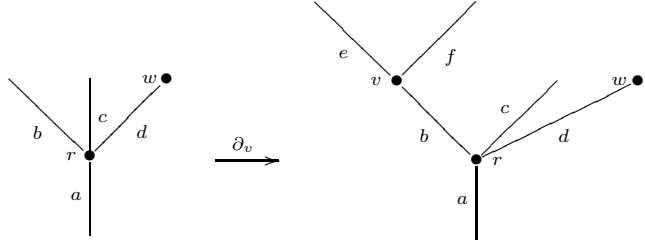
\begin{array}{ccc}
\xymatrix{\\*{\,}\ar@{-}[dr]_{b} & *{\,}\ar@{-}[d]^{c}\ar@{}[r]|{\,\,\,\,\,\,\,\,\,\,\,\, w} & *{\bullet}\ar@{-}[dl]^{d}\\
\,\ar@{}[r]|{\,\,\,\,\,\,\,\,\,\,\,\, r} & *{\bullet}\ar@{-}[d]_{a}\\
 & *{\,}}
 & \xymatrix{\\\\\ar[r]^{\partial_{v}} & *{}}
 & \xymatrix{*{\,}\ar@{-}[dr]_{e} &  & *{\,}\ar@{-}[dl]^{f}\\
\,\ar@{}[r]|{\,\,\,\,\,\,\,\,\,\,\, v} & *{\bullet}\ar@{-}[dr]_{b} &  & *{\,}\ar@{-}[dl]_{c}\ar@{}[r]|{\,\,\,\,\,\,\,\,\,\,\,\, w} & *{\bullet}\ar@{-}[dll]^{d}\\
 &  & *{\bullet}\ar@{-}[d]_{a} & \,\ar@{}[l]|{r\,\,\,\,\,\,\,\,\,\,\,}\\
 &  & *{\,}}
\end{array}\]
}

Moreover, for any tree $T$ with exactly one vertex $v$, each edge
$e$ of $T$ (necessarily outer), there is an \emph{outer face map}\[
e:\eta\rightarrow T\]
sending the unique edge of $\eta$ to $e$. 

(iii) Given a tree $T$ and an inner edge $e$ in $T$, one can obtain
a new tree $T/e$ by contracting the edge $e$. There is a canonical
map of operads $\partial_{e}:\Omega(T/e)\rightarrow\Omega(T)$ which
sends the new vertex in $T/e$ (obtained by merging the two vertices
attached to $e$) into the appropriate composition of these two vertices
in $\Omega(T)$. An arrow $\partial_{e}:T/e\rightarrow T$ in $\Omega$
of this kind is called an \emph{inner face (map).} For example

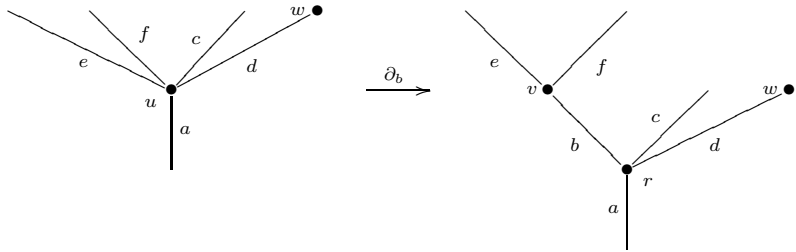
\[
\begin{array}{ccc}
\xymatrix{*{\,}\ar@{-}[rrd]_{e} & *{\,}\ar@{-}[rd]^{f} &  & *{\,}\ar@{-}[dl]_{c}\ar@{}[r]|{\,\,\,\,\,\,\,\,\, w} & *{\bullet}\ar@{-}[lld]^{d}\\
 & \,\ar@{}[r]_{\,\,\,\,\,\,\,\,\,\,\, u} & *{\bullet}\ar@{-}[d]^{a}\\
 &  & *{\,}}
 & \xymatrix{\\\ar[r]^{\partial_{b}} & *{}}
 & \xymatrix{*{\,}\ar@{-}[dr]_{e} &  & *{\,}\ar@{-}[dl]^{f}\\
\,\ar@{}[r]|{\,\,\,\,\,\,\,\,\,\,\,\,\,\, v} & *{\bullet}\ar@{-}[dr]_{b} &  & *{\,}\ar@{-}[dl]_{c}\ar@{}[r]|{\,\,\,\,\,\,\,\,\,\,\,\, w} & *{\bullet}\ar@{-}[dll]^{d}\\
 &  & *{\bullet}\ar@{-}[d]_{a} & \,\ar@{}[l]^{r\,\,\,\,\,\,\,\,\,\,\,}\\
 &  & *{\,}}
\end{array}\]

(iv) Given two trees $T$ and $T'$, any isomorphism $T\rightarrow T'$
of trees, sending inputs to inputs and output to output, of course
defines an isomorphism of operads $\Omega(T)\rightarrow\Omega(T')$,
and hence is an isomorphism $\xymatrix{T\ar[r]^{\cong} & T'}
$ in $\Omega$. For example, if $C_{n}$ denotes the corolla with just
one vertex, $n$ inputs, and one output, then we might name its input
edges $e_{1},\cdots,e_{n}$ 
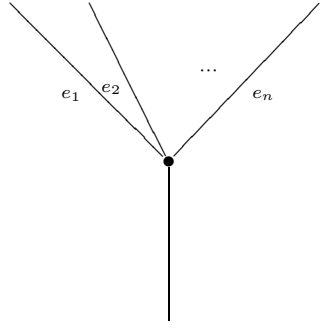
\[
\xymatrix{*{}\ar@{-}[ddrr]_{e_{1}} & *{}\ar@{-}[ddr]_{e_{2}} &  &  & *{}\ar@{-}[ddll]^{e_{n}}\\
 &  & \ar@{}[r]^{\cdots} & *{}\\
 &  & *{\bullet}\ar@{-}[dd]\\
 &  & *{}\\
 &  & *{}}
\]
 Any permutation $\varphi\in\Sigma_{n}$ defines an automorphism of
$C_{n}$ in $\Omega$. 

For a tree $T$, let its \emph{degree} $|T|$ be the number of vertices
in $T$. Then degeneracy maps decrease degree by 1, face maps (outer
or inner) increase degree by 1, and isomorphisms preserve degree.
Any map $\xymatrix{T\ar[r]^{f} & T'}
$ in $\Omega$ can be written as $f=\delta\varphi\sigma$, where $\delta$
is a composition of (inner or outer) faces, $\varphi$ is an isomorphism,
and $\sigma$ is a composition of degeneracies. (This composition
is unique up to isomorphism).

\section{Skeletal filtration}

As for any presheaf category, any dendroidal set $X$ is a colimit
of representables, of the form\[
X=\varinjlim\Omega[T]\]
 (see \cite{CWM}). We wish to refine this a little, in a way similar
to the skeletal filtration for simplicial sets. To this end, call
a dendrex $x\in X_{T}$ of shape $T$ \emph{degenerate} if there is
a surjective map $\xymatrix{T\ar[r]^{\alpha} & T'}
$ in $\Omega$ (a composition of degeneracies) such that $x=\alpha^{*}(x')$
for some $x'\in X_{T'}$. Here $\alpha$ should not be an isomorphism
of course, so that $T'$ has strictly fewer vertices then $T$ and
$\alpha$ is a non-empty composition of degeneracies. 

Given a dendroidal set $X$ we denote by $Sk_{n}(X)$ the sub dendroidal set of $X$ generated by all non-degenerate dendrices $x\in X_{T}$ where $|T|\le n$. An arbitrary
dendroidal set $X$ is clearly the colimit (union) of the sequence \[
\begin{array}{cccc}
\quad\quad\quad & Sk_{0}(X)\subseteq Sk_{1}(X)\subseteq Sk_{2}(X)\subseteq\cdots & \quad\quad\quad & (1)\end{array}\]
We call this the skeletal filtration of $X$. This filtration extends the skeletal filtration for simplicial sets
in the precise sense that for any dendroidal set $X$ and any simplicial
set $S$, there are canonical isomorphisms \[
i^{*}Sk_{n}(X)=Sk_{n}(i^{*}X)\]
and\[
i_{!}Sk_{n}(S)=Sk_{n}(i_{!}S).\]

Consider now the following diagram: \[
\xymatrix{*++{\coprod_{x,T}\partial\Omega[T]}\ar[r]\ar@{>->}[d] & *++{Sk_{n}(X)}\ar@{>->}[d]\\
\coprod_{x,T}\Omega[T]\ar[r] & Sk_{n+1}(X)}
\]
 where the sum is taken over isomorphism classes of pairs $(x,T)$ in the category of elements of $X$ where $T$ is a tree with
$n$ vertices and $x\in X_{T}$ is non-degenerate, and $\partial\Omega[T]$
is the boundary of $\Omega[T]$, i.e., the union of its faces. We call the skeletal filtration of $X$ \emph{normal} if this square is a pushout for each $n>0$. 

Following Cisinski \cite{Cisinski} we call a dendroidal set \emph{normal} if for each non-degenerate dendrex $x\in X_{T}$, the only isomorphism fixing $x$ is the identity. Cisinski (loc. cit.) proves that the normal dendroidal sets are precisely those whose skeletal filtrations are normal. 

\begin{example}
If $X$ is a simplicial set then $i_{!}(X)$ admits a normal skeletal filtration and in fact that skeletal filtration is isomorphic to the usual skeletal filtration of $X$. If $\mathcal{P}$ is a $\Sigma$-free operad then $N_d(\mathcal{P})$ is normal. In particular if $\mathcal{P}$ is the symmetrization of a planar operad then $N_d(\mathcal{P})$ is normal. 
\end{example}

\section{Inner Kan complexes}

We begin by introducing inner horns. For a tree $T$, each face map
$\partial:T'\rightarrow T$ defines a monomorphism $\Omega[T']\rightarrow\Omega[T]$
between (representable) dendroidal sets. The union (pushout) of these
subobjects is the boundary of $\Omega[T]$, denoted \[
\xymatrix{*++{\partial\Omega[T]}\ar@{>->}[r] & \Omega[T],}
\]
as above. If $e$ is an inner edge of $T$, then the union of all
the faces \emph{except} \[
\xymatrix{*+++{\partial_{e}:T/e}\ar@{>->}[r] & T}
\]
 defines a subobject of the boundary, denoted \[
\xymatrix{*++{\Lambda^{e}[T]}\ar@{>->}[r] & \Omega[T]}
,\]
 and called the \emph{inner horn} associated to $e$ (and to $T$).
This terminology and notation extends the one \[
\xymatrix{*++{\Lambda^{k}[n]}\ar@{>->}[r] & *++{\partial\Delta[n]}\ar@{>->}[r] & \Delta[n]}
\]
for simplicial sets, in the sense that \[
i_{!}(\Lambda^{k}[n])=\Lambda^{k}[i[n]]\]
\[
i_{!}(\partial\Delta[n])=\partial\Omega[i[n]]\]
as subobjects of $i_{!}(\Delta[n])=\Omega[i[n]].$

A dendroidal set $K$ is said to be a (dendroidal) \emph{inner Kan
complex} if, for any tree $T$ and any inner edge $e$ in $T$, the
map\[
K_{T}=Hom(\Omega[T],K)\rightarrow Hom(\Lambda^{e}[T],K)\]
is a surjection of sets. It is called a \emph{strict} inner Kan complex
if this map is a bijection (for any $T$ and $e$ as above). For example,
we will see (Proposition 5.3 below) that the dendroidal nerve of an
operad is always a strict inner Kan complex. This terminology is analogous
to the one introduced by Boardman and Vogt, who say a simplicial set
$X$ satisfies the restricted Kan condition if, for any $0<k<n$,
the map $Hom(\Delta[n],X)\rightarrow Hom(\Lambda^{k}[n],X)$ is a
surjection (\cite{BoardmanVogt} Definition 4.8, page 102). In more
recent work (\cite{JoyalPaper,JoyalBook}) Joyal develops the general
theory of simplicial sets satisfying the restricted Kan condition.
Joyal uses the terminology quasi-categories for such simplicial sets
so as to stress the analogy with category theory. In fact a simplicial
set $X$ is a quasi-category iff $i_{!}(X)$ is a dendroidal inner
Kan complex, and for any dendroidal inner Kan complex $K$, the restriction
$i^{*}(K)$ is a quasi-category in the sense of Joyal. 

Let us call a map $u:U\rightarrow V$ of dendroidal sets an \emph{anodyne
extension} if it can be obtained from the set of inner horn inclusions
by coproducts, pushouts, compositions, and retracts (cf, \cite{GZ},
p. 60). Then obviously, the surjectivity property for inner Kan complexes
extends to anodyne extensions, in the sense that the map of sets \[
u^{*}:Hom(V,K)\rightarrow Hom(U,K),\]
given by composition with $u$, is again surjective. Similarly, the
map $u^{*}$ is a bijection for any strict inner Kan complex. 

For a tree $T$ let $I(T)$ be the set of inner edges of $T$. For
a non-empty subset $A\subseteq I(T)$ let $\Lambda^{A}[T]$ be the
union of all faces of $\Omega[T]$ except those obtained by contracting
an edge from $A$. Note that if $A=\{ e\}$ then $\Lambda^{A}[T]=\Lambda^{e}[T]$. 

\begin{lem}
For any non-empty $A\subseteq I(T)$ the inclusion $\Lambda^{A}[T]\rightarrow\Omega[T]$
is anodyne.
\end{lem}
\begin{proof}
By induction on $n=|A|$. If $n=1$ then $\Lambda^{A}[T]\rightarrow\Omega[T]$
is an inner horn inclusion, thus anodyne. Assume the proposition holds
for $n<k$ and suppose $|A|=k$. Choose an arbitrary $e\in A$ and
put $B=A\backslash\{ e\}$. The map $\Lambda^{A}[T]\rightarrow\Omega[T]$
factors as\[
\xymatrix{\Lambda^{A}[T]\ar[r]\ar[rd] & \Lambda^{B}[T]\ar[d]\\
 & \Omega[T]}
\]
 The vertical map is anodyne by the induction hypothesis and it therefore
suffices to prove that $\Lambda^{A}[T]\rightarrow\Lambda^{B}[T]$
is anodyne. The following diagram expresses that map as a pushout\[
\xymatrix{\Lambda^{B}[T/e]\ar[r]\ar[d] & \Lambda^{A}[T]\ar[d]\\
\Omega[T/e]\ar[r] & \Lambda^{B}[T]}
\]
 and since the map $\Lambda^{B}[T/e]\rightarrow\Omega[T/e]$ is anodyne
(by the induction hypothesis), the proof is complete. 
\end{proof}
We denote by $\Lambda^{I}[T]$ the dendroidal set $\Lambda^{A}[T]$
where $A=I(T)$, that is $\Lambda^{I}[T]$ is the union of all outer
faces of $\Omega[T]$. By the above proposition the inclusion $\Lambda^{I}[T]\rightarrow\Omega[T]$
is anodyne. 

We now consider grafting of trees. For two trees $T$ and $S$, and
a leaf $l$ of $T$, let $T\circ_{l}S$ be the tree obtained by grafting
$S$ onto $T$ by identifying $l$ with the root (output edge) of
$S$. Then there are obvious inclusions $\Omega[S]\rightarrow\Omega[T\circ_{l}S]$
and $\Omega[T]\rightarrow\Omega[T\circ_{l}S]$, the pushout (union)
of which we denote by $\Omega[T]\cup_{l}\Omega[S]\rightarrow\Omega[T\circ_{l}S]$.

\begin{lem}
(Grafting) For any two trees $T$ and $S$ and any leaf $l$ of $T$,
the inclusion $\Omega[T]\cup_{l}\Omega[S]\rightarrow\Omega[T\circ_{l}S]$
is anodyne.
\end{lem}
\begin{proof}
Let us write $R=T\circ_{l}S$. The case where $T=\eta$ or $S=\eta$
is trivial, we therefore assume that this is not the case. We proceed
by induction on $n=|T|+|S|$, the sum of the degrees of $T$ and $S$.
The cases $n=0$ or $n=1$ are taken care of by our assumption that
$T\ne\eta\ne S$. For the case $n=2$ the same assumption implies
that the inclusion $\Omega[T]\cup_{l}\Omega[S]\rightarrow\Omega[R]$
is an inner horn inclusion. In any case it is anodyne. Assume then
that the result holds for $2\le n<k$ and suppose $|T|+|S|=k$. 

Recall that $\Lambda^{I}[R]$ is the union of all the outer faces
of $\Omega[R]$. First notice that $\Omega[T]\cup_{l}\Omega[S]\rightarrow\Omega[R]$
factors as\[
\xymatrix{\Omega[T]\cup_{l}\Omega[S]\ar[r]\ar[dr] & \Lambda^{I}[R]\ar[d]\\
 & \Omega[R]}
\]
 and the vertical arrow is anodyne by a previous result. We now show
that \[
\Omega[T]\cup_{l}\Omega[S]\rightarrow\Lambda^{I}[R]\]
 is anodyne by exhibiting it as a pushout of an anodyne extension.
Recall (\cite{DenSet}) that an external cluster is a vertex $v$
with the property that one of the edges adjacent to it is inner while
all the other edges adjacent to it are outer. Let $Cl(T)$ (resp.
$Cl(S)$) be the set of all external clusters in $T$ (resp. $S$)
which do not contain $l$ (resp. the root of $S$). For each $C\in Cl(T)$
the face of $\Omega[R]$ corresponding to $C$ is isomorphic to $\Omega[(T/C)\circ_{l}S]$
and the map $\Omega[T/C]\cup_{l}\Omega[S]\rightarrow\Omega[(T/C)\circ_{l}S${]}
is anodyne by the induction hypothesis. Similarly for every $C\in Cl(S)$
the face of $\Omega[R]$ that corresponds to $C$ is isomorphic to
$\Omega[T\circ_{l}(S/C)]$ and the map $\Omega[T]\cup_{l}\Omega[S/C]\rightarrow\Omega[T\circ_{l}(S/C)]$
is anodyne by the induction hypothesis. The following diagram is a
pushout\[
\xymatrix{\coprod_{C\in Cl(T)}(\Omega[T/C]\cup_{l}\Omega[S])\amalg\coprod_{C\in Cl(S)}(\Omega[T]\cup_{l}\Omega[S/C])\ar[r]\ar[d] & \Omega[T]\cup_{l}\Omega[S]\ar[d]\\
\coprod_{C\in Cl(T)}(\Omega[(T/C)\circ_{l}S])\amalg\coprod_{C\in Cl(S)}(\Omega[T\circ_{l}(S/C)])\ar[r] & \Lambda^{I}[R]}
\]
 where the map on the left is the coproduct of all of the anodyne
extensions just mentioned. Since anodyne extensions are closed under
coproducts, it follows that the map on the left of the pushout is
anodyne and thus also the one on the right, which is what we set out
to prove. This concludes the proof.  
\end{proof}
We end this section with two remarks on strict inner Kan complexes. 

\begin{prop}
The dendroidal nerve of any operad is a strict inner Kan complex.
\end{prop}
\begin{proof}
Let $\mathcal{P}$ be an operad. A dendrex $x\in N_{d}(\mathcal{P})_{T}$
is a map $x:\Omega[T]\rightarrow N_{d}(\mathcal{P})$ which is a map
of operads $\Omega(T)\rightarrow P$. If we choose a planar representative
for $T$ then $\Omega(T)$ is specifically given in terms of generators
and is a free operad. It follows that $x$ is equivalent to a labeling
of the (planar representative) $T$ as follows. The edges are labeled
by colours of $\mathcal{P}$ and the vertices are coloured by operations
in $\mathcal{P}$ where the input of such an operation is the tuple
of labels of the incoming edges to the vertex and the output is the
label of the outgoing edge from the vertex. Any inner horn $\Lambda^{e}[T]\rightarrow N_{d}(\mathcal{P})$
is easily seen to be equivalent to such a labeling of the tree $T$
and thus determines a unique filler. 
\end{proof}
\begin{prop}
Any strict inner Kan complex is $2$-coskeletal.
\end{prop}
\begin{proof}
Let $X$ be a strict inner Kan complex. Let $Y$ be any dendroidal
set and assume $f:Sk_{2}Y\rightarrow Sk_{2}X$ is given. We first
show that $f$ can be extended to a dendroidal map $\hat{f}:Y\rightarrow X$.
Suppose $f$ was extended to a map $f_{k}:Sk_{k}Y\rightarrow Sk_{k}X$
for $k\ge2$. Let $y\in Sk_{k+1}(Y)$ be a non-degenerate dendrex
and assume $y\notin Sk_{k}(Y)$. So $y\in Y_{T}$ and $T$ has exactly
$k+1$ vertices. Choose an inner horn $\Lambda^{\alpha}[T]$ (such
an inner horn exist since $k\ge2$). The set $\{\beta^{*}y\}_{\beta\ne\alpha}$
where $\beta$ runs over all faces of $T$, defines a horn $\Lambda^{\alpha}[T]\rightarrow Y$.
Since this horn factors through the $k$-skeleton of $Y$ we obtain,
by applying $f_{k}$, a horn $\Lambda^{\alpha}[T]\rightarrow X$ in
$X$ given by $\{ f(\beta^{*}y)\}_{\beta\ne\alpha}$. Let $f_{k+1}(y)\in X_{T}$
be the unique filler of that horn. By construction we have for each
$\beta\ne\alpha$ that \[
\beta^{*}f_{k+1}(y)=f(\beta^{*}y)\]
 it thus remains to show the same for $\alpha$. The dendrices $f(\alpha^{*}y)$
and $\alpha^{*}f_{k+1}(y)$ both have the same boundary and they are
both of shape $S$ where $S$ has $k$ vertices. Since $k\ge2$, $S$
has an inner face, but then it follows that both $f(\alpha^{*}y)$
and $\alpha^{*}f_{k+1}(y)$ are fillers for the same inner horn in
$X$ which proves that they are equal. By repeating the process for
all dendrices in $Sk_{k+1}(Y)$ it follows that $f_{k}$ can be extended
to $f_{k+1}:Sk_{k+1}(Y)\rightarrow Sk_{k+1}(X)$. This holds for all
$k\ge2$ which implies that $f$ can be extended to $\hat{f}:Y\rightarrow X$.
To show uniqueness of $\hat{f}$ assume that $g$ is another extension
of $f$. Suppose it has been shown that $\hat{f}$ and $g$ agree
on all dendrices of shape $T$ where $T$ has at most $k$ vertices,
and let $y\in X_{S}$ be a dendrex of shape $S$ where $S$ has $k+1$
vertices. But then the dendrices $\hat{f}(y)$ and $g(y)$ are dendrices
in $X$ that have the same boundary. Since $k\ge2$ it follows that
these dendrices are both fillers for the same inner horn and so are
the same. This proves that $\hat{f}=g$. 
\end{proof}

\section{The operad generated by an inner Kan complex}

We recall that $\tau_{d}:dSet\rightarrow Operad$ denotes the left
adjoint to the dendroidal nerve functor $N_{d}$. In this section,
we will give a more explicit description of the operad $\tau_{d}(X)$
in the case where $X$ is an inner Kan complex. This description extends
the one in \cite{BoardmanVogt} of the category generated by a simplicial
set satisfying the restricted Kan condition. It will lead to a proof
of the following converse of Proposition 5.3.

\begin{thm}
For any strict inner Kan complex $X$, the canonical map $X\rightarrow N_{d}(\tau_{d}(X))$
is an isomorphism.
\end{thm}
Proposition 5.3 and Theorem 6.1 together state that a dendroidal set
is a strict inner Kan complex iff it is the nerve of an operad. 

Consider an inner Kan complex $X$. For the description of $\tau_{d}(X)$,
we first fix some notation. For each $n\ge0$ let $C_{n}$ be the
$n$-corolla:
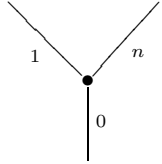
\[
\xymatrix{*{}\ar@{-}[dr]_{1} &  & *{}\ar@{-}[dl]^{n}\\
 & *{\bullet}\ar@{-}[d]^{0}\\
 & *{}}
\]
and for each $0\le i\le n$ recall that $i:\eta\rightarrow C_{n}$
denotes the obvious (outer face) map in $\Omega$ that sends the unique
edge of $\eta$ to the edge $i$ in $C_{n}$. An element $f\in X_{C_{n}}$
will be denoted by 
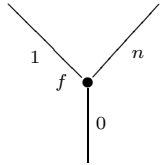
\[
\xymatrix{*{}\ar@{-}[dr]_{1} &  & *{}\ar@{-}[dl]^{n}\\
\ar@{}[r]|{\quad f} & *{\bullet}\ar@{-}[d]^{0}\\
 & *{}}
\]
If $C_{n}^{'}$ is another $n$-corolla together with an isomorphism
$\alpha:C_{n}^{'}\rightarrow C_{n}$ then we will usually write $f$
again instead of $\alpha^{*}(f)$. We will use this convention quite
often in the coming definitions and constructions, and in each case
there will be an obvious choice for the isomorphism $\alpha$ given
by the planar representation of the trees at question, which will
usually not be mentioned. 

\begin{defn}
Let $X$ be an inner Kan complex and let $f,g\in X_{C_{n}}$, $n\ge0$.
For $1\le i\le n$ we say that \emph{$f$ is homotopic to $g$ along
the edge $i$,} and write $f\sim_{i}g$, if there is a dendrex $H$
of shape 
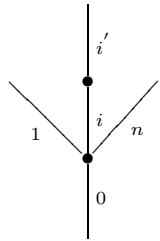
\[
\xymatrix{ & *{}\ar@{-}[d]^{i^{'}}\\
*{}\ar@{-}[dr]_{1} & *{\bullet}\ar@{-}[d]^{i} & *{}\ar@{-}[dl]^{n}\\
 & *{\bullet}\ar@{-}[d]^{0}\\
 & *{}}
\]
whose three faces are:\[
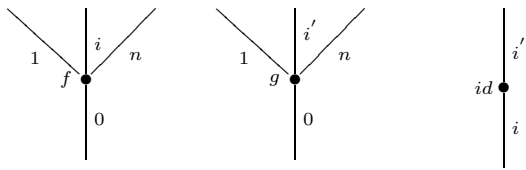
\begin{array}{ccc}
\xymatrix{*{}\ar@{-}[dr]_{1} & *{}\ar@{-}[d]^{i} & *{}\ar@{-}[dl]^{n}\\
\ar@{}[r]|{\quad\,\,\, f} & *{\bullet}\ar@{-}[d]^{0}\\
 & *{}}
 & \quad\xymatrix{*{}\ar@{-}[dr]_{1} & *{}\ar@{-}[d]\ar@{}[d(0.6)]^{i^{'}} & *{}\ar@{-}[dl]^{n}\\
\ar@{}[r]|{\quad\,\,\, g} & *{\bullet}\ar@{-}[d]^{0}\\
 & *{}}
 & \quad\xymatrix{ & *{}\ar@{-}[d]^{i^{'}}\\
\ar@{}[r]|{\quad\,\,\, id} & *{\bullet}\ar@{-}[d]^{i}\\
 & *{}}
\end{array}\]
where the third one denotes a degeneracy. Similarly we will say that
$f$ is homotopic to $g$ along the edge $0$ and write $f\sim_{0}g$
if there is a dendrex of shape
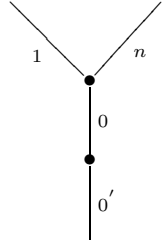
\[
\xymatrix{*{}\ar@{-}[dr]_{1} &  & *{}\ar@{-}[dl]^{n}\\
 & *{\bullet}\ar@{-}[d]^{0}\\
 & *{\bullet}\ar@{-}[d]^{0^{'}}\\
 & *{}}
\]
whose three faces are:\[
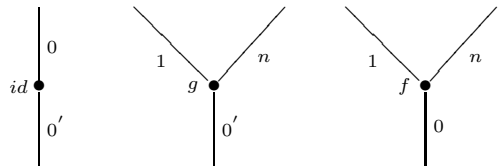
\begin{array}{ccc}
\xymatrix{ & *{}\ar@{-}[d]^{0}\\
\ar@{}[r]|{\quad\,\,\, id} & *{\bullet}\ar@{-}[d]^{0^{'}}\\
 & *{}}
\quad & \xymatrix{*{}\ar@{-}[dr]_{1} &  & *{}\ar@{-}[dl]^{n}\\
\ar@{}[r]|{\quad\,\,\, g} & *{\bullet}\ar@{-}[d]^{0^{'}}\\
 & *{}}
\quad & \xymatrix{*{}\ar@{-}[dr]_{1} &  & *{}\ar@{-}[dl]^{n}\\
\ar@{}[r]|{\quad\,\,\, f} & *{\bullet}\ar@{-}[d]^{0}\\
 & *{}}
\end{array}\]
When $f\sim_{i}g$ for some $0\le i\le n$ we will refer to the corresponding
$H$ as a \emph{homotopy} from $f$ to $g$ along $i$ and will sometimes
write $H:f\sim_{i}g$. 
\end{defn}
\begin{prop}
Let $X$ be an inner Kan complex. For each $0\le i\le n$ the relation
$\sim_{i}$ on the set $X_{C_{n}}$ is an equivalence relation.
\end{prop}
\begin{proof}
First we prove reflexivity. For $1\le i\le n$ let \[
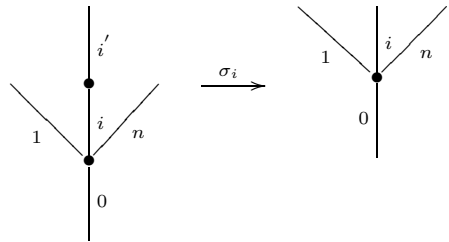
\begin{array}{ccc}
\xymatrix{ & *{}\ar@{-}[d]^{i^{'}}\\
*{}\ar@{-}[dr]_{1} & *{\bullet}\ar@{-}[d]^{i} & *{}\ar@{-}[dl]^{n}\\
 & *{\bullet}\ar@{-}[d]^{0}\\
 & *{}}
 & \xymatrix{\\\ar[r]^{\sigma_{i}} & *{}\\
}
 & \xymatrix{*{}\ar@{-}[dr]_{1} & *{}\ar@{-}[d]^{i} & *{}\ar@{-}[dl]^{n}\\
 & *{\bullet}\ar@{-}[d]_{0}\\
 & *{}}
\end{array}\]
and for $i=0$ let\[
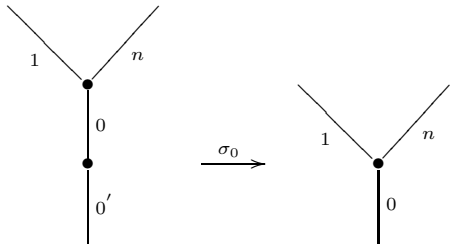
\begin{array}{ccc}
\xymatrix{*{}\ar@{-}[dr]_{1} &  & *{}\ar@{-}[dl]^{n}\\
 & *{\bullet}\ar@{-}[d]^{0}\\
 & *{\bullet}\ar@{-}[d]^{0^{'}}\\
 & *{}}
 & \xymatrix{\\\\\ar[r]^{\sigma_{0}} & *{}\\
}
 & \xymatrix{\\*{}\ar@{-}[dr]_{1} &  & *{}\ar@{-}[dl]^{n}\\
 & *{\bullet}\ar@{-}[d]^{0}\\
 & *{}}
\end{array}\]
be the obvious degeneracies. It then follows that for any $f\in X_{C_{n}}$
the dendrex $\sigma_{i}^{*}(f)$ is a homotopy from $f$ to $f$,
thus $f\sim_{i}f$. 

To prove symmetry assume $f\sim_{i}g$ for some $1\le i\le n$ and
let $H_{fg}$ be a homotopy from $f$ to $g$ along $i$. Consider
the tree $T$:
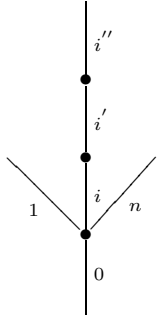
\[
\xymatrix{\\ & *{}\ar@{-}[d]^{i^{''}}\\
 & *{\bullet}\ar@{-}[d]^{i^{'}}\\
*{}\ar@{-}[dr]_{1} & *{\bullet}\ar@{-}[d]^{i} & *{}\ar@{-}[dl]^{n}\\
 & *{\bullet}\ar@{-}[d]^{0}\\
 & *{}}
\]
 For the inner horn $\Lambda^{i}[T]$, corresponding to the edge $i$
in the tree above, we now describe a map $\Lambda^{i}[T]\rightarrow X$.
Such a map is given by specifying three dendrices in $X$ of certain
shapes such that their faces match in a suitable way. We describe
this map by explicitly writing the mentioned dendrices and their faces:\\
\\
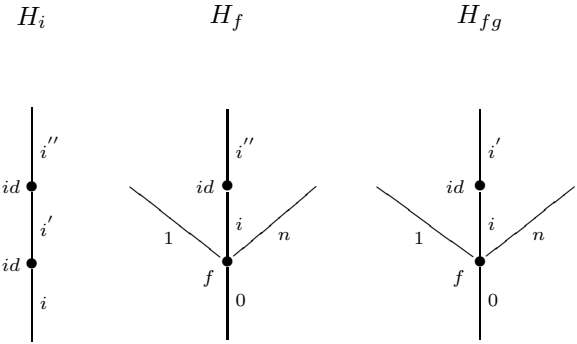
\[
\begin{array}{ccc}
\xymatrix{ & H_{i}\\
 & *{}\ar@{-}[d]^{i^{''}}\\
\ar@{}[r]|{\quad\quad id}\ar@{}[d] & *{\bullet}\ar@{-}[d]^{i^{'}}\\
*{}\ar@{}[r]|{\quad\quad id} & *{\bullet}\ar@{-}[d]^{i}\\
 & *{}\\
}
 & \quad\xymatrix{ & H_{f}\\
 & *{}\ar@{-}[d]^{i^{''}}\\
*{}\ar@{-}[dr]_{1}\ar@{}[r]|{\quad\quad id} & *{\bullet}\ar@{-}[d]^{i} & *{}\ar@{-}[dl]^{n}\\
*{}\ar@{}[r]_{\quad\,\,\quad f} & *{\bullet}\ar@{-}[d]^{0}\\
 & *{}\\
}
\quad & \xymatrix{ & H_{fg}\\
 & *{}\ar@{-}[d]^{i^{'}}\\
*{}\ar@{-}[dr]_{1}\ar@{}[r]|{\quad\quad id} & *{\bullet}\ar@{-}[d]^{i} & *{}\ar@{-}[dl]^{n}\\
*{}\ar@{}[r]_{\,\,\quad\quad f} & *{\bullet}\ar@{-}[d]^{0}\\
 & *{}\\
}
\end{array}\]
with inner faces of these dendrices:\[
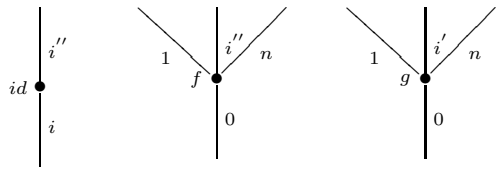
\begin{array}{ccc}
\xymatrix{ & *{}\ar@{-}[d]^{i^{''}}\\
\ar@{}[r]|{\quad\,\, id} & *{\bullet}\ar@{-}[d]^{i}\\
 & *{}}
\quad & \xymatrix{*{}\ar@{-}[dr]_{1} & *{}\ar@{-}[d]^{i^{''}} & *{}\ar@{-}[dl]^{n}\\
\ar@{}[r]|{\quad\,\,\, f} & *{\bullet}\ar@{-}[d]^{0}\\
 & *{}}
\quad & \xymatrix{*{}\ar@{-}[dr]_{1} & *{}\ar@{-}[d]^{i^{'}} & *{}\ar@{-}[dl]^{n}\\
\ar@{}[r]|{\quad\,\,\, g} & *{\bullet}\ar@{-}[d]^{0}\\
 & *{}}
\end{array}\]
where $H_{i}$ is a double degeneracy of $i$, $H_{f}$ is a homotopy
from $f$ to $f$ (along the branch $i$) and $H_{fg}$ is the given
homotopy from $f$ to $g$. It is easily checked that the faces indeed
match so that we have a horn in $X$. Let $x$ be a filler for that
horn and consider $H_{gf}=\partial_{i}^{*}(x)$. This dendrex can
be pictured as
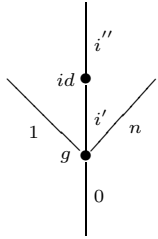
\[
\xymatrix{ & *{}\ar@{-}[d]^{i^{''}}\\
*{}\ar@{-}[dr]_{1}\ar@{}[r]|{\quad\,\,\, id} & *{\bullet}\ar@{-}[d]^{i'} & *{}\ar@{-}[dl]^{n}\\
\ar@{}[r]|{\quad\,\,\, g} & *{\bullet}\ar@{-}[d]^{0}\\
 & *{}}
\]
with inner face:
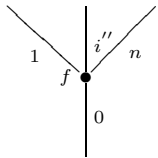
\[
\xymatrix{*{}\ar@{-}[dr]_{1} & *{}\ar@{-}[d]^{i^{''}} & *{}\ar@{-}[dl]^{n}\\
\ar@{}[r]|{\quad\,\,\, f} & *{\bullet}\ar@{-}[d]^{0}\\
 & *{}}
\]
and is thus a homotopy from $g$ to $f$ along $i$, so that $g\sim_{i}f$.
For $i=0$ a similar proof works. 

To prove transitivity let $f\sim_{i}g$ and $g\sim_{i}h$ for $1\le i\le n$.
Let $H_{fg}$ be a homotopy from $f$ to $g$ and let $H_{gh}$ be
a homotopy from $g$ to $h$. We again consider the tree $T$:
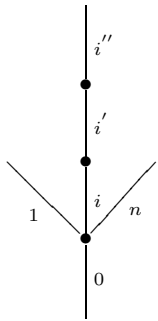
\[
\xymatrix{\\ & *{}\ar@{-}[d]^{i^{''}}\\
 & *{\bullet}\ar@{-}[d]^{i^{'}}\\
*{}\ar@{-}[dr]_{1} & *{\bullet}\ar@{-}[d]^{i} & *{}\ar@{-}[dl]^{n}\\
 & *{\bullet}\ar@{-}[d]^{0}\\
 & *{}}
\]
This time we look at $\Lambda^{i^{'}}[T]$. The following describes
a map $\Lambda^{i^{'}}[T]\rightarrow X$ in $X$:\[
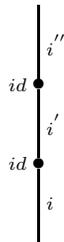
\begin{array}{ccc}
\xymatrix{ & H_{i}\\
 & *{}\ar@{-}[d]^{i^{''}}\\
\ar@{}[r]|{\quad\quad id} & *{\bullet}\ar@{-}[d]^{i^{'}}\\
\ar@{}[r]|{\quad\quad id} & *{\bullet}\ar@{-}[d]^{i}\\
 & *{}}
\quad & 
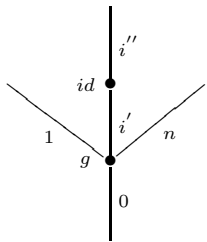
\xymatrix{ & H_{gh}\\
 & *{}\ar@{-}[d]^{i^{''}}\\
*{}\ar@{-}[dr]_{1}\ar@{}[r]|{\quad\quad id} & *{\bullet}\ar@{-}[d]^{i^{'}} & *{}\ar@{-}[dl]^{n}\\
\ar@{}[r]|{\quad\quad g} & *{\bullet}\ar@{-}[d]^{0}\\
 & *{}}
\quad & 
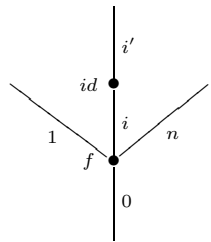
\xymatrix{ & H_{fg}\\
 & *{}\ar@{-}[d]^{i'}\\
*{}\ar@{-}[dr]_{1}\ar@{}[r]|{\quad\quad id} & *{\bullet}\ar@{-}[d]^{i} & *{}\ar@{-}[dl]^{n}\\
\ar@{}[r]|{\quad\quad f} & *{\bullet}\ar@{-}[d]^{0}\\
 & *{}}
\end{array}\]
with inner faces being:\[
\xyC{34pt}\begin{array}{ccc}
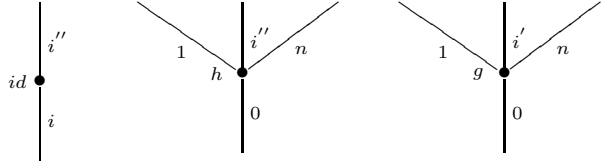
\xymatrix{ & *{}\ar@{-}[d]^{i^{''}}\\
\ar@{}[r]|{\,\,\quad\quad id} & *{\bullet}\ar@{-}[d]^{i}\\
 & *{}}
\quad & \xymatrix{*{}\ar@{-}[dr]_{1} & *{}\ar@{-}[d]^{i^{''}} & *{}\ar@{-}[dl]^{n}\\
\ar@{}[r]|{\quad\quad h} & *{\bullet}\ar@{-}[d]^{0}\\
 & *{}}
\quad & \xymatrix{*{}\ar@{-}[dr]_{1} & *{}\ar@{-}[d]^{i^{'}} & *{}\ar@{-}[dl]^{n}\\
\ar@{}[r]|{\quad\quad g} & *{\bullet}\ar@{-}[d]^{0}\\
 & *{}}
\end{array}\]
Let $x$ be a filler for that horn and let $H_{fh}=\partial_{i^{'}}^{*}(x)$.
This dendrex can be pictured as follows:
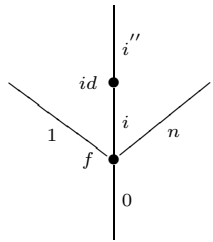
\[
\xymatrix{ & *{}\ar@{-}[d]^{i^{''}}\\
*{}\ar@{-}[dr]_{1}\ar@{}[r]|{\quad\quad id} & *{\bullet}\ar@{-}[d]^{i} & *{}\ar@{-}[dl]^{n}\\
\ar@{}[r]|{\quad\quad f} & *{\bullet}\ar@{-}[d]^{0}\\
 & *{}}
\]
with inner face:
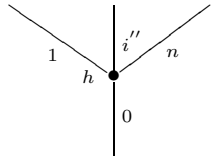
\[
\xymatrix{*{}\ar@{-}[dr]_{1} & *{}\ar@{-}[d]^{i^{''}} & *{}\ar@{-}[dl]^{n}\\
\ar@{}[r]|{\quad\quad h} & *{\bullet}\ar@{-}[d]^{0}\\
 & *{}}
\]
and is thus a homotopy from $f$ to $h$ so that $f\sim_{i}h$. The
proof for $i=0$ is similar.
\end{proof}
\begin{lem}
Let $X$ be an inner Kan complex. The relations $\sim_{0},\cdots,\sim_{n}$
on $X_{C_{n}}$ are all equal.
\end{lem}
\begin{rem}
On the basis of this lemma, we will later just write $f\sim g$ instead
of $f\sim_{i}g$.
\end{rem}
\begin{proof}
Suppose $H:f\sim_{i}g$ for $1\le i\le n$ and let $1\le i<j\le n$.
We consider the tree T:\[
\xyC{15pt}
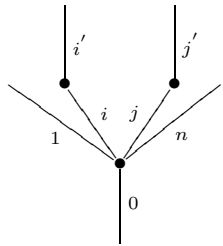
\xymatrix{ & *{}\ar@{-}[d]^{i^{'}} &  & *{}\ar@{-}[d]^{j^{'}}\\
*{}\ar@{-}[drr]_{1} & *{\bullet}\ar@{-}[dr]^{i} &  & *{\bullet}\ar@{-}[dl]_{j} & *{}\ar@{-}[dll]^{n}\\
 &  & *{\bullet}\ar@{-}[d]^{0}\\
 &  & *{}}
\]
and the inner horn $\Lambda^{i}[T]$. The following then describes
a map $\Lambda^{i}[T]\rightarrow X$:\\
\\
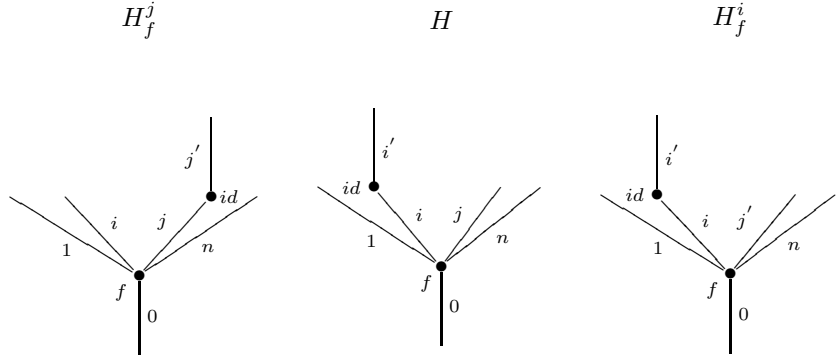
\[
\begin{array}{ccc}
\xyC{15pt}\xymatrix{ &  & H_{f}^{j}\\
 &  &  & *{}\ar@{-}[d]_{j^{'}}\\
*{}\ar@{-}[drr]_{1} & *{}\ar@{-}[dr]^{i} &  & *{\bullet}\ar@{-}[dl]_{j} & *{}\ar@{-}[dll]^{n}\ar@{}[l]|{id\,\,\,}\\
 & \ar@{}[r]_{\quad\,\,\, f} & *{\bullet}\ar@{-}[d]^{0}\\
 &  & *{}}
\quad & \xymatrix{ &  & H\\
 & *{}\ar@{-}[d]^{i^{'}}\\
*{}\ar@{-}[drr]_{1}\ar@{}[r]|{\,\,\, id} & *{\bullet}\ar@{-}[dr]^{i} &  & *{}\ar@{-}[dl]_{j} & *{}\ar@{-}[dll]^{n}\\
 & \ar@{}[r]_{\,\,\,\quad f} & *{\bullet}\ar@{-}[d]^{0}\\
 &  & *{}}
\quad & \xymatrix{ &  & H_{f}^{i}\\
 & *{}\ar@{-}[d]^{i^{'}}\\
*{}\ar@{-}[drr]_{1}\ar@{}[r]|{\,\,\, id} & *{\bullet}\ar@{-}[dr]^{i} &  & *{}\ar@{-}[dl]_{j^{'}} & *{}\ar@{-}[dll]^{n}\\
 & \ar@{}[r]_{\quad\,\,\, f} & *{\bullet}\ar@{-}[d]^{0}\\
 &  & *{}}
\end{array}\]
where $H_{f}^{j}:f\sim_{j}f$ and $H_{f}^{i}:f\sim_{i}f$. The inner
faces of the three dendrices are\[
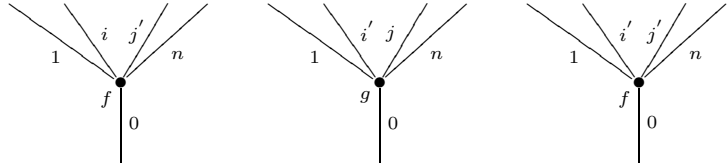
\begin{array}{ccc}
\xyC{15pt}\xymatrix{*{}\ar@{-}[drr]_{1} & *{}\ar@{-}[dr]^{i} & \ar@{}[d(0.7)]^{j^{'}} & *{}\ar@{-}[dl] & *{}\ar@{-}[dll]^{n}\\
 & \ar@{}[r]_{\quad f} & *{\bullet}\ar@{-}[d]^{0}\\
 &  & *{}}
 & \quad\xymatrix{*{}\ar@{-}[drr]_{1} & *{}\ar@{-}[dr]^{i^{'}} &  & *{}\ar@{-}[dl]_{j} & *{}\ar@{-}[dll]^{n}\\
 & \ar@{}[r]_{\quad g} & *{\bullet}\ar@{-}[d]^{0}\\
 &  & *{}}
\quad & \xymatrix{*{}\ar@{-}[drr]_{1} & *{}\ar@{-}[dr]^{i^{'}} & \ar@{}[d(0.7)]^{j^{'}} & *{}\ar@{-}[dl] & *{}\ar@{-}[dll]^{n}\\
 & \ar@{}[r]_{\quad f} & *{\bullet}\ar@{-}[d]^{0}\\
 &  & *{}}
\end{array}\]
 Let $x$ be a filler for this horn, then $\partial_{i}^{*}(x)$ is
the following dendrex
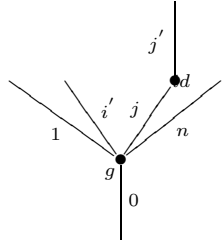
\[
\xymatrix{ &  &  & *{}\ar@{-}[d]_{j^{'}}\\
*{}\ar@{-}[drr]_{1} & *{}\ar@{-}[dr]^{i^{'}} &  & *{\bullet}\ar@{-}[dl]_{j} & *{}\ar@{-}[dll]^{n}\ar@{}[l]|{id\,\quad\,}\\
 & \ar@{}[r]_{\quad\,\, g} & *{\bullet}\ar@{-}[d]^{0}\\
 &  & *{}}
\]
with inner face:
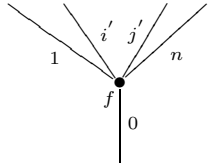
\[
\xymatrix{*{}\ar@{-}[drr]_{1} & *{}\ar@{-}[dr]^{i^{'}} & \ar@{}[d(0.7)]^{j^{'}} & *{}\ar@{-}[dl] & *{}\ar@{-}[dll]^{n}\\
 & \ar@{}[r]_{\quad\,\, f} & *{\bullet}\ar@{-}[d]^{0}\\
 &  & *{}}
\]
and is thus a homotopy from $g$ to $f$ along the $j$-th branch.
Thus $g\sim_{j}f$ and so $f\sim_{j}g$ as well. The other cases to
be considered follow in a similar way.
\end{proof}
Given an inner Kan complex $X$ and vertices $x_{1},\cdots,x_{n},x\in X_{\eta}$,
let us write \[
X(x_{1},\cdots,x_{n};x)\subseteq X(C_{n})\]
for the set of dendrices $x$ of shape $C_{n}$ with $0^{*}(x)=x$
and $i^{*}(x)=x_{i}$ for $i=0,\cdots,n$. Here $i:\eta\rightarrow C_{n}$
denotes the map in $\Omega$ sending the unique edge of $\eta$ to
the one of $C_{n}$ with name $i$. The equivalence relation $\sim$
on $X(C_{n})$ given by the preceding lemma defines a quotient of
$X(C_{n})$ which we will denote by \[
Ho(X)(x_{1},\cdots,x_{n};x)=X(x_{1},\cdots,x_{n};x)/\sim.\]
 This defines a coloured collection $Ho(X)$, and a canonical quotient
map of collections $Sk_{1}(X)\rightarrow Ho(X)$. We will now proceed
to prove the following.

\begin{prop}
There is a unique structure of a (symmetric, coloured) operad on $Ho(X)$
for which the map of collections $Sk_{1}(X)\rightarrow Ho(X)$ extends
to a map of dendroidal sets $X\rightarrow N_{d}(Ho(X))$. The latter
map is an isomorphism whenever $X$ is a \emph{strict} inner Kan \emph{}complex\emph{. }
\end{prop}
To prepare for the proof of this proposition, we begin by defining
the composition operations $\circ_{i}$ of the operad $Ho(X)$. Let
$X$ be an inner Kan complex and let $f\in X_{C_{n}}$ and $g\in X_{C_{m}}$
be two dendrices in $X$. We will say that a dendrex $h\in X_{C_{n+m-1}}$
is a $\circ_{i}$-composition of $f$ and $g$ if there is a dendrex
$\gamma$ in $X$ as follows:\[
\xyC{25pt}\xyR{25pt}
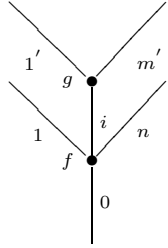
\xymatrix{*{}\ar@{-}[dr]_{1^{'}} &  & *{}\ar@{-}[dl]^{m^{'}}\\
*{}\ar@{-}[dr]_{1}\ar@{}[r]|{\,\,\quad g} & *{\bullet}\ar@{-}[d]^{i} & *{}\ar@{-}[dl]^{n}\\
\ar@{}[r]|{\,\,\quad f} & *{\bullet}\ar@{-}[d]^{0}\\
 & *{}}
\]
with inner face
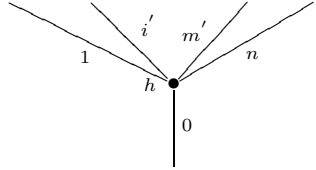
\[
\xymatrix{*{}\ar@{-}[drr]_{1} & *{}\ar@{-}[dr]^{i^{'}} & \ar@{}[d(0.7)]^{m^{'}\,} & *{}\ar@{-}[dl] & *{}\ar@{-}[dll]^{n}\\
 & \ar@{}[r]|{\quad\,\, h} & *{\bullet}\ar@{-}[d]^{0}\\
 &  & *{}}
\]
We will denote this situation by $h\sim f\circ_{i}g$ and call $\gamma$
a \emph{witness} for the composition. 

\begin{rem}
Notice that for $1\le i\le n$ we have by definition that $H:f\sim_{i}g$
iff $H$ is a witness for the composition $g\sim f\circ_{i}id$. Similarly
for $i=0$ we have that $H:f\sim_{0}g$ iff $H$ is a witness for
the composition $g\sim id\circ f$. 
\end{rem}
\begin{lem}
In an inner Kan complex $X$, if $h\sim f\circ_{i}g$ and $h'\sim f\circ_{i}g$
then $h\sim h'$. 
\end{lem}
\begin{proof}
Let $\gamma$ be a witness for the composition $h\sim f\circ_{i}g$
and $\gamma'$ one for the composition $h'\sim f\circ_{i}g$. We consider
the tree $T$:
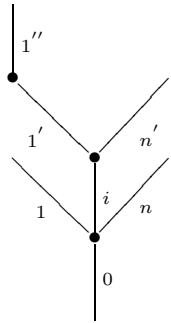
\[
\xymatrix{*{}\ar@{-}[d]^{1^{''}}\\
*{\bullet}\ar@{-}[dr]_{1^{'}} &  & *{}\ar@{-}[dl]^{n^{'}}\\
*{}\ar@{-}[dr]_{1} & *{\bullet}\ar@{-}[d]^{i} & *{}\ar@{-}[dl]^{n}\\
 & *{\bullet}\ar@{-}[d]^{0}\\
 & *{}}
\]
and the inner horn $\Lambda^{i}[T]$. Let $H_{g}:g\sim_{i}g$ and
consider the following map $\Lambda^{i}[T]\rightarrow X$\[
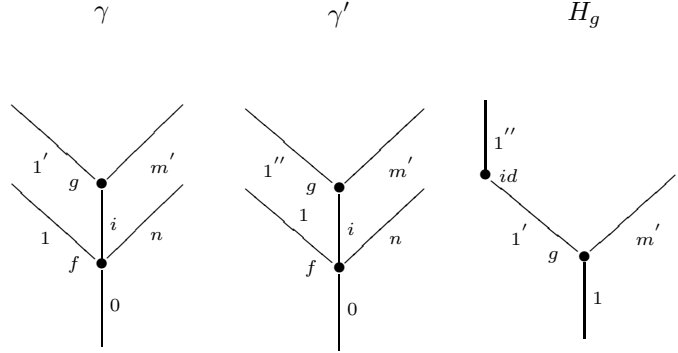
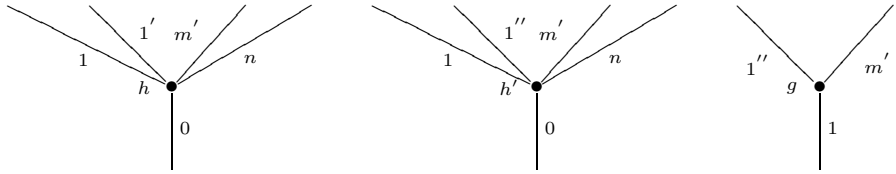
\begin{array}{ccc}
\xymatrix{ & \gamma\\
*{}\ar@{-}[dr]_{1^{'}} &  & *{}\ar@{-}[dl]^{m^{'}}\\
*{}\ar@{-}[dr]_{1}\ar@{}[r]|{\quad\,\, g} & *{\bullet}\ar@{-}[d]^{i} & *{}\ar@{-}[dl]^{n}\\
\ar@{}[r]|{\quad\,\, f} & *{\bullet}\ar@{-}[d]^{0}\\
 & *{}}
\quad & \xymatrix{ & \gamma'\\
*{}\ar@{-}[dr]_{1^{''}} &  & *{}\ar@{-}[dl]^{m^{'}}\\
*{}\ar@{-}[dr]^{1}\ar@{}[r]|{\quad\,\,\, g} & *{\bullet}\ar@{-}[d]^{i} & *{}\ar@{-}[dl]^{n}\\
\ar@{}[r]|{\quad\,\,\, f} & *{\bullet}\ar@{-}[d]^{0}\\
 & *{}}
\quad & \xymatrix{ & H_{g}\\
*{}\ar@{-}[d]^{1^{''}}\\
*{\bullet}\ar@{-}[dr]_{1^{'}} & \ar@{}[l]|{id\quad\quad} & *{}\ar@{-}[dl]^{m^{'}}\\
\ar@{}[r]|{\quad\,\,\, g} & *{\bullet}\ar@{-}[d]^{1}\\
 & *{}}
\end{array}\]
with inner faces\[
\begin{array}{ccc}
\xymatrix{*{}\ar@{-}[drr]_{1} & *{}\ar@{-}[dr]^{1^{'}} &  & *{}\ar@{-}[dl]_{m^{'}} & *{}\ar@{-}[dll]^{n}\\
 & \ar@{}[r]|{\quad h} & *{\bullet}\ar@{-}[d]^{0}\\
 &  & *{}}
\quad & \xymatrix{*{}\ar@{-}[drr]_{1} & *{}\ar@{-}[dr]^{1^{''}} &  & *{}\ar@{-}[dl]_{m^{'}} & *{}\ar@{-}[dll]^{n}\\
 & \ar@{}[r]|{\quad h'} & *{\bullet}\ar@{-}[d]^{0}\\
 &  & *{}}
\quad & \xymatrix{*{}\ar@{-}[dr]_{1^{''}} &  & *{}\ar@{-}[dl]^{m^{'}}\\
\ar@{}[r]|{\quad g} & *{\bullet}\ar@{-}[d]^{1}\\
 & *{}}
\end{array}\]
Let $x$ be a filler for this horn. The face $\partial_{i}^{*}x$
is then the dendrex
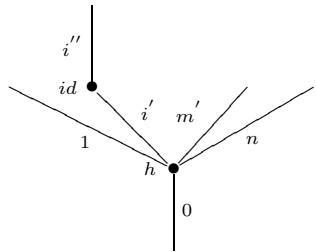
\[
\xymatrix{ & *{}\ar@{-}[d]_{i^{''}}\\
*{}\ar@{-}[drr]_{1}\ar@{}[r]|{\,\,\quad id} & *{\bullet}\ar@{-}[dr]^{i^{'}} &  & *{}\ar@{-}[dl]_{m^{'}} & *{}\ar@{-}[dll]^{n}\\
 & \ar@{}[r]|{\,\,\quad h} & *{\bullet}\ar@{-}[d]^{0}\\
 &  & *{}}
\]
whose inner face is 
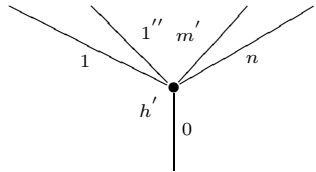
\[
\xymatrix{*{}\ar@{-}[drr]_{1} & *{}\ar@{-}[dr]^{1^{''}} &  & *{}\ar@{-}[dl]_{m^{'}} & *{}\ar@{-}[dll]^{n}\\
 & \ar@{}[r]_{\quad\,\, h^{'}} & *{\bullet}\ar@{-}[d]^{0}\\
 &  & *{}}
\]
which proves that $h\sim h'$.
\end{proof}
\begin{lem}
In an inner Kan complex $X$, let $f\sim f'$ and $g\sim g'$. If
$h\sim f\circ_{i}g$ and $h'\sim f'\circ_{i}g'$ then $h\sim h'$. 
\end{lem}
\begin{proof}
Let $H$ be a homotopy from $f$ to $f'$ along the edge $i$, $H'$
a homotopy from $g'$ to $g$ along the root, and $\gamma$ a witness
for the composition $h\sim f\circ_{i}g$. We now consider the tree
T:
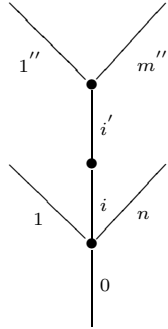
\[
\xymatrix{*{}\ar@{-}[dr]_{1^{''}} &  & *{}\ar@{-}[dl]^{m^{''}}\\
 & *{\bullet}\ar@{-}[d]^{i^{'}}\\
*{}\ar@{-}[dr]_{1} & *{\bullet}\ar@{-}[d]^{i} & *{}\ar@{-}[dl]^{n}\\
 & *{\bullet}\ar@{-}[d]^{0}\\
 & *{}}
\]
and the inner horn $\Lambda^{i}[T]$. The following is then a map
$\Lambda^{i}[T]\rightarrow X$ in $X$:\\
\\
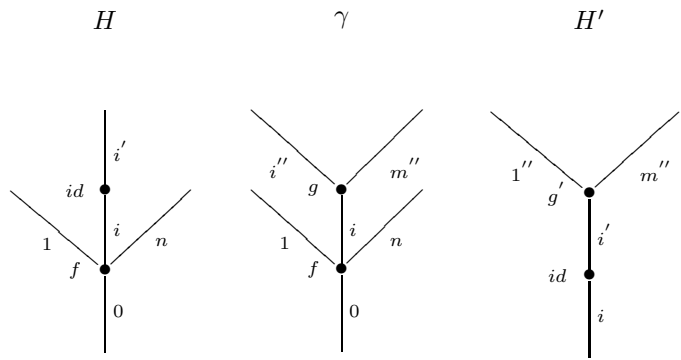
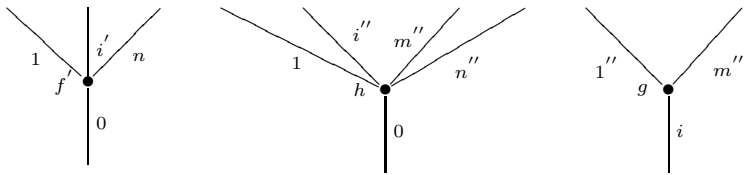
\[
\xyC{25pt}\xyR{25pt}\begin{array}{ccc}
\xymatrix{ & H\\
 & *{}\ar@{-}[d]^{i^{'}}\\
*{}\ar@{-}[dr]_{1}\ar@{}[r]|{\,\,\quad id} & *{\bullet}\ar@{-}[d]^{i} & *{}\ar@{-}[dl]^{n}\\
\ar@{}[r]|{\,\,\quad f} & *{\bullet}\ar@{-}[d]^{0}\\
 & *{}}
\quad & \xymatrix{ & \gamma\\
*{}\ar@{-}[dr]_{i^{''}} &  & *{}\ar@{-}[dl]^{m^{''}}\\
*{}\ar@{-}[dr]_{1}\ar@{}[r]|{\,\,\quad g} & *{\bullet}\ar@{-}[d]^{i} & *{}\ar@{-}[dl]^{n}\\
\ar@{}[r]|{\,\,\quad f} & *{\bullet}\ar@{-}[d]^{0}\\
 & *{}}
\quad & \xymatrix{ & H'\\
*{}\ar@{-}[dr]_{1^{''}} &  & *{}\ar@{-}[dl]^{m^{''}}\\
\ar@{}[r]|{\,\,\quad g^{'}} & *{\bullet}\ar@{-}[d]^{i^{'}}\\
\ar@{}[r]|{\,\,\quad id} & *{\bullet}\ar@{-}[d]^{i}\\
 & *{}}
\end{array}\]
with inner faces:\[
\xyC{25pt}\xyR{25pt}\begin{array}{ccc}
\xymatrix{*{}\ar@{-}[dr]_{1} & *{}\ar@{-}[d]^{i^{'}} & *{}\ar@{-}[dl]^{n}\\
\ar@{}[r]|{\quad\,\, f^{'}} & *{\bullet}\ar@{-}[d]^{0}\\
 & *{}}
\quad & \xymatrix{*{}\ar@{-}[drr]_{1} & *{}\ar@{-}[dr]^{i^{''}} & \ar@{}[d(0.7)]^{m^{''}} & *{}\ar@{-}[dl] & *{}\ar@{-}[dll]^{n^{''}}\\
 & \ar@{}[r]|{\quad\, h} & *{\bullet}\ar@{-}[d]^{0}\\
 &  & *{}}
\quad & \xymatrix{*{}\ar@{-}[dr]_{1^{''}} &  & *{}\ar@{-}[dl]^{m^{''}}\\
\ar@{}[r]|{\quad\, g} & *{\bullet}\ar@{-}[d]^{i}\\
 & *{}}
\end{array}\]
The missing face of a filler for this horn is then:\[
\xyC{25pt}\xyR{25pt}
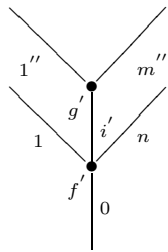
\xymatrix{*{}\ar@{-}[dr]_{1^{''}} &  & *{}\ar@{-}[dl]^{m^{''}}\\
*{}\ar@{-}[dr]_{1}\ar@{}[r]_{\quad\quad g^{'}} & *{\bullet}\ar@{-}[d]^{i^{'}} & *{}\ar@{-}[dl]^{n}\\
\ar@{}[r]_{\quad\quad f^{'}} & *{\bullet}\ar@{-}[d]^{0}\\
 & *{}}
\]
 with inner face
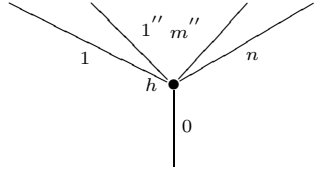
\[
\xymatrix{*{}\ar@{-}[drr]_{1} & *{}\ar@{-}[dr]^{1^{''}} &  & *{}\ar@{-}[dl]_{m^{''}} & *{}\ar@{-}[dll]^{n}\\
 & \ar@{}[r]|{\quad\,\,\, h} & *{\bullet}\ar@{-}[d]^{0}\\
 &  & *{}}
\]
which proves that $h\sim f'\circ_{i}g'$, and thus by the previous
result also that $h\sim h'$. 
\end{proof}
We now proceed to prove Proposition 5.6:

\begin{proof}
(of Proposition 5.6) Lemma 5.8 implies that for \[
[f]\in Ho(X)(x_{1},\cdots,x_{n};x)\]
 and \[
[g]\in Ho(X)(y_{1},\cdots,y_{m};x_{i})\]
 the assignment \[
[f]\circ_{i}[g]=[f\circ_{i}g]\]
 is well-defined. This provides the $\circ_{i}$ operations in the
operad $Ho(X)$. The $\Sigma_{n}$ actions are defines as follows.
Given a permutation $\sigma\in\Sigma_{n}$ let $\sigma:C_{n}\rightarrow C_{n}$
be the obvious induced map in $\Omega$. The map $\sigma^{*}:X_{C_{n}}\rightarrow X_{C_{n}}$
restricts to a function \[
\sigma^{*}:X(x_{1},\cdots,x_{n};x)\rightarrow X(x_{\sigma(1)},\cdots,x_{\sigma(n)};x)\]
 and it is trivial to verify that this map respects the homotopy relation.
We thus obtain a map \[
\sigma^{*}Ho(X)(x_{1},\cdots,x_{n};x)\rightarrow Ho(X)(x_{\sigma(1)},\cdots,x_{\sigma(n)};x).\]
We now need to show that these structure maps make the coloured collection
$Ho(X)$ into an operad. The verification is simple and we exemplify
it by proving associativity. Let $[f]\in Ho(X)(x_{1},\cdots,x_{n};x)$,
$[g]\in Ho(X)(y_{1},\cdots,y_{m};x_{i})$ and $[h]\in Ho(X)(z_{1},\cdots,z_{k};y_{m})$.
We need to prove that $[f]\circ([g]\circ[h])=([f]\circ[g])\circ[h]$
(with the obvious choice for subscripts on the $\circ$) which is
the same as showing that $f\circ(g\circ h)\sim(f\circ g)\circ h$
for any choice for compositions $\psi\sim g\circ h$ and $\varphi\sim f\circ g$.
Consider the tree $T$ given by
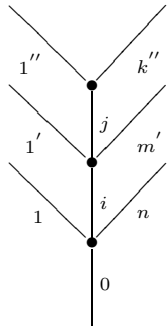
\[
\xymatrix{*{}\ar@{-}[dr]_{1^{''}} &  & *{}\ar@{-}[dl]^{k^{''}}\\
*{}\ar@{-}[dr]_{1^{'}} & *{\bullet}\ar@{-}[d]^{j} & *{}\ar@{-}[dl]^{m^{'}}\\
*{}\ar@{-}[dr]_{1} & *{\bullet}\ar@{-}[d]^{i} & *{}\ar@{-}[dl]^{n}\\
 & *{\bullet}\ar@{-}[d]^{0}\\
 & *{}}
\]
and consider the anodyne extension $\Lambda^{I}[T]\rightarrow\Omega[T]$,
cf. Lemma 5.1. The two given compositions $\psi\sim g\circ h$ and
$\varphi\sim f\circ g$ define a map $\Lambda^{I}[T]\rightarrow X$
depicted by\[
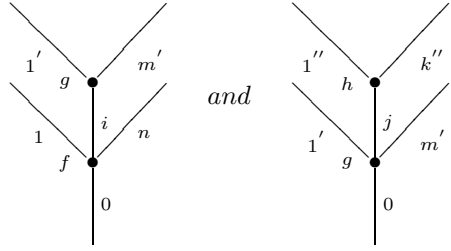
\begin{array}{ccc}
\xymatrix{*{}\ar@{-}[dr]_{1^{'}} &  & *{}\ar@{-}[dl]^{m^{'}}\\
*{}\ar@{-}[dr]_{1}\ar@{}[r]|{\quad g} & *{\bullet}\ar@{-}[d]^{i} & *{}\ar@{-}[dl]^{n}\\
\ar@{}[r]|{\quad f} & *{\bullet}\ar@{-}[d]^{0}\\
 & *{}}
 & \xymatrix{\\and\\
}
 & \xymatrix{*{}\ar@{-}[dr]_{1^{''}} &  & *{}\ar@{-}[dl]^{k^{''}}\\
*{}\ar@{-}[dr]_{1^{'}}\ar@{}[r]|{\quad h} & *{\bullet}\ar@{-}[d]^{j} & *{}\ar@{-}[dl]^{m^{'}}\\
\ar@{}[r]|{\quad g} & *{\bullet}\ar@{-}[d]^{0}\\
 & *{}}
\end{array}\]
whose inner faces are respectively $\psi$ and $\varphi$. Let $x\in X_{T}$
be a dendrex extending this map. Let $c:C_{m}\rightarrow T$ be the
map obtained by contracting both $i$ and $j$ and $\rho=c^{*}x$.
It now follows that $\partial_{i}^{*}x$ is a witness for the the
composition $\rho\sim\psi\circ h$ and $\partial_{j}^{*}x$ is a witness
for the composition $\rho\sim f\circ\varphi$. That proves associativity.
The other axioms for an operad follow in a similar manner. 

Next, let us show that the quotient map $q:Sk_{1}(X)\rightarrow Ho(X)$
extends to a map $q:X\rightarrow N_{d}(Ho(X))$ of dendroidal sets.
Since we already know that $N_{d}(X)$ is $2$-coskeletal, it suffices
to give its values for dendrices $x\in X_{T}$ where $T$ is a tree
with two vertices. Let $e$ be the inner edge of this tree. Then $\xymatrix{*++{\Lambda^{e}[T]}\ar@{>->}[r] & \Omega[T]\ar[r]^{x} & X}
$ factors through $Sk_{1}(X)$, so its composition $\Lambda^{e}[T]\rightarrow N_{d}(Ho(X))$
with $q$ has a unique extension (Proposition 5.3), which we take
to be $q(x):\Omega[T]\rightarrow N_{d}(Ho(X))$. This defines $q:Sk_{2}(X)\rightarrow Sk_{2}(N_{d}(Ho(X))),$
and hence all of $q:X\rightarrow N_{d}(Ho(X))$ by $2$-coskeletality,
as said. 

Finally, when $X$ is itself a strict inner Kan complex, then the
homotopy relation is the identity relation, so $Sk_{1}(X)\rightarrow Ho(X)$
is the identity map. Since $X$ and $N_{d}(Ho(X))$ are now both strict
inner Kan complexes, the extension $q:X\rightarrow N_{d}(Ho(X))$
is an isomorphism. 
\end{proof}
The following Proposition, together with Proposition 6.6, now provide
the proof of Theorem 6.1.

\begin{prop}
For any inner Kan complex $X$, the natural map $\tau_{d}(X)\rightarrow Ho(X)$
is an isomorphism of operads. 
\end{prop}
\begin{proof}
It suffices to prove that the map $q:X\rightarrow N_{d}(Ho(X))$ of
Proposition 6.6 has the universal property of the unit of the adjunction.
This means that for any operad $\mathcal{P}$ and any map $\varphi:X\rightarrow N_{d}(\mathcal{P})$,
there is a unique map of operads $\psi:Ho(X)\rightarrow\mathcal{P}$
for which $N_{d}(\psi)q=\varphi$. But $\varphi$ induces a map $Ho(X)\rightarrow Ho(N_{d}(\mathcal{P}))$
for which\[
\xymatrix{Sk_{1}(X)\ar[r]^{\varphi}\ar[d]^{q_{X}} & Sk_{1}N_{d}(\mathcal{P})\ar[d]^{q_{\mathcal{P}}}\\
Ho(X)\ar[r]^{Ho(\varphi)} & Ho(N_{d}(\mathcal{P}))}
\]
commutes, and $Ho(N_{d}(\mathcal{P}))=\mathcal{P}$ while $q_{\mathcal{P}}$
is the identity as we have seen in (the proof of) Proposition 6.6.
So $Ho(\varphi)$ in fact defines a map $\psi:Ho(X)\rightarrow\mathcal{P}$
of collections. It is easily seen that $\psi$ is a map of operads.
It is unique because $q_{X}$ is surjective. 
\end{proof}

\section{Homotopy coherent nerves of operads}

In this section, we assume $\mathcal{E}$ is a monoidal model category
with a cofibrant unit $I$. We also assume that $\mathcal{E}$ is
equipped with an \emph{interval} in the sense of \cite{GenWCon}.
Such an interval is given by maps \[
\xymatrix{I\ar@<2pt>[r]^{0}\ar@<-2pt>[r]_{1} & H\ar[r]^{\epsilon} & I}
\]
and \[
\xymatrix{H\otimes H\ar[r]^{\,\,\,\,\,\,\vee} & H}
\]
 satisfying certain conditions. In particular, $H$ is an interval
in Quillen's sense (\cite{homotopicalAlg}), so $0$ and $1$ together
define a cofibration $I\coprod I\rightarrow H$, and $\epsilon$ is
a weak equivalence. Such an interval $H$ allows one to construct
for each (coloured) operad $\mathcal{P}$ in $\mathcal{E}$ a ''Boardman-Vogt''
resolution $W_{H}(\mathcal{P})\rightarrow\mathcal{P}$. Each operad
in $Set$ can be viewed as an operad in $\mathcal{E}$ (via the functor
$Set\rightarrow\mathcal{E}$ which preserves coproducts and sends
the one-point set to $I$), and hence has such a Boardman-Vogt resolution.
When we apply this to the operads $\Omega(T)$, we obtain the \emph{homotopy
coherent dendroidal nerve $hcN_{d}(\mathcal{P})$} of any operad $\mathcal{P}$
in $\mathcal{E}$, as the dendroidal set given by \[
hcN_{d}(\mathcal{P})_{T}=Hom(W_{H}(\Omega(T)),\mathcal{P})\]
where the Hom is that of operads in $\mathcal{E}$. See \cite{DenSet}
for a more detailed description and examples. Our goal here is to
prove the following result.

\begin{thm}
Let $\mathcal{P}$ be an operad in $\mathcal{E}$, with the property
that for each sequence $c_{1},\cdots,c_{n};c$ of colours of $\mathcal{P}$,
the object $\mathcal{P}(c_{1},\cdots,c_{n};c)$ is fibrant. Then $hcN_{d}(\mathcal{P})$
is an inner Kan complex. 
\end{thm}
\begin{rem}
As explained in \cite{DenSet}, our construction of the dendroidal
homotopy coherent nerve specializes to that of the homotopy coherent
nerve of an $\mathcal{E}$-enriched category, and for the case where
$\mathcal{E}$ is the category of topological spaces or simplicial
sets, one recovers the classical definition (\cite{CordierPorter}). In particular,
as a special case of Theorem 7.1, one obtains that for an $\mathcal{E}$-enriched
category with fibrant Hom objects (in other words, for a locally fibrant
$\mathcal{E}$-enriched category), its homotopy coherent nerve is
a quasi-category in the sense of Joyal. This result was proved, for
the case where $\mathcal{E}$ is simplicial sets, by Cordier and Porter
in \cite{CordierPorter}. 
\end{rem}
Before embarking on the proof of Theorem 7.1, we need to be a bit
more explicit about the operads of the form $W_{H}\Omega(T)$ involved
in the definition of the homotopy coherent nerve. Recall first of
all the functor\[
Symm:Operad(\mathcal{E})_{\pi}\rightarrow Operad(\mathcal{E})\]
 which is left adjoint to the forgetful functor from symmetric operads
to non-symmetric (i.e., planar) ones. If $T$ is an object in $\Omega$
and $\bar{T}$ is a chosen planar representative of $T$, then $\bar{T}$
naturally describes a planar operad $\Omega(\bar{T})$ for which $\Omega(T)=Symm(\Omega(\bar{T}))$.
Since the $W$-construction commutes with symmetrization (as one readily
verifies), it follows that \[
W_{H}(T)=Symm(W_{H}\Omega(\bar{T})).\]

This latter operad $W_{H}\Omega(\bar{T})$ is easily described explicitly.
The colours of \[
W_{H}(\Omega(\bar{T}))\]
 are the colours of $\Omega(\bar{T})$, i.e., the edges of $T$. By
a \emph{signature,} we mean a sequence $e_{1},\cdots,e_{n};e_{0}$
of edges. Given a signature $\sigma=(e_{1},\cdots,e_{n};e_{0})$,
we have that $W_{H}(\Omega(\bar{T}))(\sigma)=0$ whenever $\Omega(\bar{T})(\sigma)=\phi$.
And if $\Omega(\bar{T})(\sigma)\ne\phi$, there is a subtree $T_{\sigma}$
of $T$ (and a corresponding planar subtree $\bar{T}_{\sigma}$ of
$\bar{T}$) whose leaves are $e_{1},\cdots,e_{n},$ and whose root
is $e_{0}$. Then \[
W_{H}\Omega(\bar{T})(e_{1},\cdots,e_{n};e_{0})=\bigotimes_{f\in i(\sigma)}H,\]
where $i(\sigma)$ is the set of \emph{inner} edges of $T_{\sigma}$
(or of $\bar{T}_{\sigma}$). (This last tensor product is to be thought
of as the ''space'' of assignments of lengths to inner edges in
$\bar{T}_{\sigma}$; it is the unit if $i(\sigma)$ is empty.)

\begin{rem}
The composition operations in the operad $W_{H}\Omega(\bar{T})$ are
given in terms of the $\circ_{i}$-operations as follows. For signatures
$\sigma=(e_{1},\cdots,e_{n};e_{0})$ and $\rho=(f_{1},\cdots,f_{m};e_{i})$,
the composition map \[
\begin{array}{cccc}
\quad & \xymatrix{\Omega(\bar{T})(e_{1},\cdots,e_{n};e_{0})\otimes\Omega(\bar{T})(f_{1},\cdots,f_{m};e_{i})\ar[d]^{\circ_{i}} & \ar@{}[d]|{(1)}\\
\Omega(\bar{T})(e_{1},\cdots,e_{i-1},f_{1},\cdots,f_{m},e_{i+1},\cdots,e_{n};e_{0}) & *{}}
 & \quad\end{array}\]
is the following one. The trees $\bar{T}_{\sigma}$ and $\bar{T}_{\rho}$
can be grafted along $e_{i}$ to form $\bar{T}_{\sigma}\circ_{e_{i}}\bar{T}_{\rho}$,
again a planar subtree of $\bar{T}$. In fact\[
\bar{T}_{\sigma}\circ_{e_{i}}\bar{T}_{\rho}=\bar{T}_{\sigma\circ_{i}\rho}\]
where $\sigma\circ_{i}\rho$ is the signature $(e_{1},\cdots,e_{i-1},f_{1},\cdots,f_{m},e_{i+1},\cdots,e_{n};e_{0})$,
and for the sets of inner edges we have \[
i(\sigma\circ_{i}\rho)=i(\sigma)\cup i(\rho)\cup\{ e_{i}\}.\]

The composition map in $(1)$ now is the map\[
\xymatrix{H^{\otimes i(\sigma)}\otimes H^{\otimes i(\rho)}\ar@{..>}[r]\ar[d]^{\cong} & H^{\otimes i(\sigma\circ_{i}\rho)}\ar[d]^{\cong}\\
H^{\otimes i(\sigma)\cup i(\rho)}\otimes I\ar[r]^{id\otimes1} & H^{\otimes i(\sigma)\cup i(\rho)}\otimes H}
\]
where $1:I\rightarrow H$ is one of the ''endpoints'' of the interval
$H$, as above. 

This description of the operad $W_{H}\Omega(\bar{T})$ is functorial
in the planar tree $T$. In particular, we note that for an inner
edge $e$ of $T$, the tree $T/e$ inherits a planar structure $\overline{T/e}$
from $\bar{T}$, and $W_{H}\Omega(\overline{T/e})\rightarrow W_{H}\Omega(\bar{T})$
is the natural map assigning length $0$ to the edge $e$ whenever
it occurs (in a subtree given by a signature). 
\end{rem}
\begin{proof}
(Of Theorem 7.1) Consider a tree $T$ and an inner edge $e$ in $T$.
We want to solve the extension problem\[
\xymatrix{*+++{\Lambda^{e}[T]}\ar[r]^{\varphi\,\,\,\,\,\,}\ar@{>->}[d] & hcN_{d}(\mathcal{P})\\
\Omega[T]\ar@{..>}[ur]}
\]

Fix a planar representative $\bar{T}$ of $T$. Then the desired map
$\psi:\Omega[T]\rightarrow hcN_{d}(\mathcal{P})$ corresponds to a
map of planar operads \[
\hat{\psi}:W_{H}\Omega(\bar{T})\rightarrow\mathcal{P}.\]
Each face $S$ of $T$ inherits a planar structure $\bar{S}$ from
$\bar{T}$, and the given map $\varphi:\Lambda^{e}[T]\rightarrow hcN_{d}(\mathcal{P})$
corresponds to a map of operads in $\mathcal{E}$,\[
\hat{\varphi}:W_{H}(\Lambda^{e}[T])\rightarrow\mathcal{P},\]
where $W_{H}(\Lambda^{e}[T])$ denotes the colimit of operads in $\mathcal{E}$,
\[
\begin{array}{ccc}
\quad\quad\quad\quad\quad\quad & W_{H}(\Lambda^{e}[T])=colim\, W(\Omega(\bar{S})) & \quad\quad\quad\quad\quad\quad(2)\end{array}\]
over all but one of the faces of $T$. In other words, $\varphi$
corresponds to a compatible family of maps \[
\hat{\varphi}_{S}:W_{H}(\Omega(\bar{S}))\rightarrow\mathcal{P}.\]
Let us now show the existence of an operad map $\hat{\psi}$ extending
the $\hat{\varphi}_{S}$ for all faces $S\ne T/e$. First, the colours
of $\Omega(\bar{T})$ are the same as those of the colimit in $(2)$,
so we already have a map $\psi_{0}=\varphi_{0}$ on colours:\[
\psi_{0}:(\textrm{Edges of }T)\rightarrow(\textrm{Colours of $\mathcal{P}$}).\]
Next, if $\sigma=(e_{1},\cdots,e_{n};e_{0})$ is a signature of $T$
for which $W_{H}(\Omega(\bar{T}))\ne\phi$, and if $T_{\sigma}\subseteq T$
is not all of $T$, then $T_{\sigma}$ is contained in an outer face
$S$ of $T$. So $W_{H}(\Omega(\bar{T}))(\sigma)=W_{H}(\Omega(\bar{T}_{\sigma}))(\sigma)=W_{H}(\Omega(\bar{S}))(\sigma)$,
and we already have a map \[
\hat{\varphi}_{S}(\sigma):W_{H}(\Omega(\bar{T}))(\sigma)\rightarrow\mathcal{P}(\sigma),\]
given by $\hat{\varphi}_{S}:W_{H}(\Omega(\bar{S}))\rightarrow\mathcal{P}$.
Thus, the only part of the operad map $\hat{\psi}:W_{H}(\Omega(\bar{T}))\rightarrow\mathcal{P}$
not determined by $\varphi$ is the one for the signature $\tau$
where $T_{\tau}=T$; i.e., $\tau=(e_{1},\cdots,e_{n};e_{0})$ where
$e_{1},\cdots,e_{n}$ are all the input edges of $\bar{T}$ (in the
planar order) and $e_{0}$ is the output edge. For this signature,
$\hat{\psi}(\tau)$ is to be a map\[
\hat{\psi}:W_{H}(\Omega(\bar{T})(\tau)=H^{\otimes i(\tau)}\rightarrow\mathcal{P}(\tau)\]
which (i) is compatible with the $\hat{\psi}(\sigma)=\hat{\varphi}_{S}(\sigma)$
for other signatures $\sigma;$ and (ii) together with these $\hat{\psi}(\sigma)$
respects operad composition. The first condition determines $\hat{\psi}(\tau)$
on the subobject of $H^{\otimes i(\tau)}$ which is given by a value
$0$ on one of the tensor-factors marked by an edge $e_{i}$ \emph{other}
than the given $e$. The second condition determines $\hat{\psi}(\tau)$
on the subobject of $H^{\otimes i(\tau)}$ which is given by a value
$1$ on one of the factors. Thus, if we write $1$ for the map $\xymatrix{*++{I}\ar@{>->}[r]^{1} & H}
$ and $\xymatrix{*++{\partial H}\ar@{>->}[r] & H}
$ for the map $I\coprod I\rightarrow H$, and define $\xymatrix{*++{\partial H^{\otimes k}}\ar@{>->}[r] & H^{\otimes k}}
$ by the Leibniz rule (i.e., $\partial(A\otimes B)=\partial(A)\otimes B\cup A\otimes\partial(b)$),
then the problem of finding $\hat{\psi}(\tau)$ comes down to an extension
problem of the form\[
\xymatrix{*+++{\partial(H^{\otimes i(\sigma)-\{ e\}}\otimes H)\cup H^{\otimes i(\sigma)-\{ e\}}\otimes I}\ar[r]\ar@{>->}[d] & \mathcal{P}(\tau)\\
H^{\otimes i(\sigma)-\{ e\}}\otimes H\ar[r]^{\cong} & H^{\otimes i(\sigma)}\ar@{..>}[u]^{\hat{\psi}(\sigma)}}
\]
This extension problem has a solution, because $\mathcal{P}(\tau)$
is fibrant by assumption, and because the left hand map is a trivial
cofibration (by repeated use of the push-out product axiom for monoidal
model categories). This concludes the proof of the theorem.
\end{proof}

\section{Grothendieck construction for dendroidal sets}

Let $\mathbb{S}$ be a Cartesian category. A functor $X:\mathbb{S}^{op}\rightarrow dSet$
is called a \emph{diagram} of dendroidal sets. In \cite{DenSet} a
construction was given of the dendroidal set $\int_{\mathbb{S}}X$.
This construction was then applied to the specific diagram of dendroidal
sets $X:Set^{op}\rightarrow dSet$, where for a set $A$, $X(A)$
was the dendroidal set of weak $n$-categories having $A$ as set
of objects. The dendroidal set $\int_{\mathbb{S}}X$ was defined to
be the dendroidal set of weak $n$-categories. Our aim in this section
is to prove that for a given diagram of dendroidal sets $X:\mathbb{S}^{op}\rightarrow dSet$,
if each $X(S)$ is an inner Kan complex then $\int_{\mathbb{S}}X$
is also an inner Kan complex. For the convenience of the reader we
repeat here the definition of $\int_{\mathbb{S}}X$. 

It will be convenient to consider dendroidal collections. A dendroidal
collection is a collection of sets $X=\{ X_{T}\}_{T\in\Omega}$. Each
dendroidal set has an obvious underlying dendroidal collection. A
map of dendroidal collections $X\rightarrow Y$ is a collection of
functions $\{ X_{T}\rightarrow Y_{T}\}_{T\in\Omega}$. Given a Cartesian
category $\mathbb{S}$, consider the dendroidal nerve $N_{d}(\mathbb{S})$
where $\mathbb{S}$ is regarded as an operad via the Cartesian structure.
There is a natural way of associating an object of $\mathbb{S}$ with
each dendrex of $N_{d}(\mathbb{S})$. For a tree $T$ in $\Omega$,
let $leaves(T)$ be the set of leaves of $T$, and for a leaf $l$,
write $l:\eta\rightarrow T$ also for the map sending the unique edge
in $\eta$ to $l$ in $T$. Then, since $\mathbb{S}$ is assumed to
have finite products, each dendrex $t\in N_{d}(\mathbb{S})_{T}$ defines
an object \[
in(t)=\prod_{l\in leaves(T)}l^{*}(t)\]
in $\mathbb{S}$. Notice that if $\alpha:S\rightarrow T$ is a composition
of face maps, then by using the canonical symmetries and the projections
in $\mathbb{S}$ there is a canonical arrow $in(\alpha):in(t)\rightarrow in(\alpha^{*}t)$
for any $t\in X_{T}$. 

\begin{defn}
Let $X:\mathbb{S}^{op}\rightarrow dSet$ be a diagram of dendroidal
sets. The dendroidal set $\int_{\mathbb{S}}X$ is defined as follows.
A dendrex $\Omega[T]\rightarrow\int_{\mathbb{S}}X$ is a pair $(t,x)$
such that $t\in N_{d}(\mathbb{S})_{T}$ and $x$ is a map of dendroidal
collections \[
x:\Omega[T]\rightarrow\coprod_{S\in ob(\mathbb{S})}X(S)\]
satisfying the following conditions. For each $r\in\Omega[T]_{R}$
(that is an arrow $r:R\rightarrow T$), we demand that $x(r)\in X(in(r^{*}t))$.
Furthermore we demand the following compatibility conditions to hold.
For any $r\in\Omega[T]_{R}$ and any map $\alpha:\xymatrix{*++{U}\ar@{>->}[r] & R}
$ in $\Omega$ \[
\alpha^{*}(x(r))=X(in(\alpha))x(\alpha^{*}(r)).\]

\end{defn}
\begin{thm}
Let $X:\mathbb{S}^{op}\rightarrow dSet$ be a diagram of dendroidal
sets. If for any $S\in ob(\mathbb{S})$ the dendroidal set $X(S)$
is a (strict) inner Kan complex then so is $\int_{\mathbb{S}}X$. 
\end{thm}
\begin{proof}
Let $T$ be a tree and $e$ an inner edge. We consider the extension
problem\[
\xymatrix{*++{\Lambda^{e}[T]}\ar[r]\ar@{>->}[d] & \int_{\mathbb{S}}X\\
\Omega[T]\ar@{..>}[ur]}
\]
The horn $\Lambda^{e}[T]\rightarrow\int_{\mathbb{S}}X$ is given by
a compatible collection $\{(r,x_{R}):\Omega[R]\rightarrow\int_{\mathbb{S}}X\}_{R\ne T/e}$.
We wish to construct a dendrex $(t,x_{T}):\Omega[T]\rightarrow\int_{\mathbb{S}}X$
extending this family. First notice that the collection $\{ r\}_{R\ne T/e}$
is an inner horn $\Lambda^{e}[T]\rightarrow N_{d}(\mathbb{S})$ (actually
this horn is obtained by composing with the obvious projection $\int_{\mathbb{S}}X\rightarrow N_{d}(\mathbb{S})$
sending a dendrex $(t,x)$ to $t$). We already know $N_{d}(\mathbb{S})$
to be an inner Kan complex (actually a strict inner Kan complex) and
thus there is a (unique) filler $t\in N_{d}(\mathbb{S})_{T}$ for
the horn $\{ r\}_{R\ne T/e}$. We now wish to define a map of dendroidal
collections $x_{T}:\Omega[T]\rightarrow\coprod_{S\in ob(\mathbb{S})}X(S)$
that will extend the given maps $x_{R}$ for $R\ne T/e$. This condition
already determines the value of $x_{T}$ for any dendrex $r:U\rightarrow T$
other then $id:T\rightarrow T$ and $\alpha:T/e\rightarrow T$, since
for each such $r$, the tree $U$ factors through one of the faces
$R\ne T/e$. To determine $x_{T}(id_{T})$ and $x_{T}(\alpha)$ consider
the family $\{ y_{R}=x_{R}(id:R\rightarrow R)\}_{R\ne T/e}$. By definition
we have that $y_{R}\in X(in(r))_{R}$. For each such $R$ let $\alpha_{R}:R\rightarrow T$
be the corresponding face map in $\Omega$. Since $\alpha^{*}t=r$
we obtain the map $in(\alpha_{R}):in(r)\rightarrow in(t)$. We can
now pull back the collection $\{ y_{R}\}_{R\ne T/e}$ using $X(in(\alpha_{R}))$
to obtain a collection $\{ z_{R}=X(in(\alpha_{R}))(y_{R})\}_{R\ne T/e}$.
This collection is now a horn $\Lambda^{e}[T]\rightarrow X(in(T))$
(this follows from the compatibility conditions in the definition
of $\int_{\mathbb{S}}X$). Since $X(in(t))$ is inner Kan there is
a filler $u\in X(in(t))_{T}$ for that horn. We now define $x_{T}(id:T\rightarrow T)=u$
and $x_{T}(\alpha:T/e\rightarrow T)=\alpha^{*}(u)$. Notice that since
$e$ is inner we have that $in(t)=in(\alpha)$ and thus the image
of these dendrices are in the correct dendroidal set, namely $X(in(t))$.
It follows from our construction that this makes $(t,x_{T})$ a dendrex
$\Omega[T]\rightarrow\int_{\mathbb{S}}X$ which extends the given
horn. This concludes the proof. 
\end{proof}

\section{The exponential property}

Our aim in this section is to prove the following theorem concerning
the closed monoidal structure of dendroidal sets. 

\begin{thm}
Let $K$ and $X$ be dendroidal sets, and assume $X$ is normal. If $K$ is a (strict) inner
Kan complex, then so is $\underline{Hom}_{dSet}(X,K)$.
\end{thm}
The internal \emph{Hom} here is defined by the universal property,
giving a bijective correspondence between maps $Y\otimes X\rightarrow K$
and $Y\rightarrow\underline{Hom}(X,K)$ for any dendroidal set $Y$,
and natural in $Y$. We recall from Section 2 that $\otimes$ is defined
in terms of the Boardman-Vogt tensor product of operads. We remind
the reader that for two (coloured) operads $\mathcal{P}$ and $\mathcal{Q}$
with respective sets of colours $C$ and $D$, this tensor product
operad $\mathcal{P}\otimes_{BV}\mathcal{Q}$ has the product $C\times D$
as its set of colours, and is described in terms of generators and
relations as follows. The operations in $\mathcal{P}\otimes_{BV}\mathcal{Q}$
are generated by the operations\[
p\otimes d\in\mathcal{P}\otimes_{BV}\mathcal{Q}((c_{1},d),\cdots,(c_{n},d);(c,d))\]
for any $p\in\mathcal{P}(c_{1},\cdots,c_{n};c)$ and any $d\in D$,
and \[
c\otimes q\in\mathcal{P}\otimes_{BV}\mathcal{Q}((c,d_{1}),\cdots,(c,d_{m});(c,d))\]
for any $q\in\mathcal{Q}(d_{1},\cdots,d_{m};d)$ and any $c\in C$.
The \emph{relations} between these state, first of all, that for fixed
$c\in C$ and $d\in D$, the maps $p\mapsto p\otimes d$ and $c\mapsto c\otimes q$
are maps of operads. Secondly, there is an \emph{interchange law}
stating that, for $p$ and $q$ as above, the composition $p\otimes d(c\otimes q,\cdots,c\otimes q)$
in\[
\mathcal{P}\otimes_{BV}\mathcal{Q}((c_{1},d_{1}),\cdots,(c_{1},d_{m}),\cdots,(c_{n},d_{1}),\cdots,(c_{n},d_{m});(c,d))\]
 and $c\otimes q(p\otimes d,\cdots,p\otimes d)$ in\[
\mathcal{P}\otimes_{BV}\mathcal{Q}((c_{1},d_{1}),\cdots,(c_{n},d_{1}),\cdots,(c_{1},d_{m}),\cdots,(c_{n},d_{m});(c,d))\]
are mapped to each other by the obvious permutation $\tau\in\Sigma_{n\times m}$
which puts the two sequences of input colours in the same order. The
tensor product of dendroidal sets is then uniquely determined (up
to isomorphism) by the fact that it preserves colimits in each variable
separately, together with the identity\[
\Omega[S]\otimes\Omega[T]=N_{d}(\Omega(S)\otimes_{BV}\Omega(T))\]
stated in Section 2, which gives the tensor product of two representable
dendroidal sets. 

First of all, let us prove that Theorem 9.1 follows by a standard
argument from the following proposition.

\begin{prop}
For any two objects $S$ and $T$ of $\Omega$, and any inner edge
$e$ in $S$, the map\[
\xymatrix{*+++{\Lambda^{e}[S]\otimes\Omega[T]\cup\Omega[S]\otimes\partial\Omega[T]}\ar@{>->}[r] & \Omega[S]\otimes\Omega[T]}
\]
is an anodyne extension. 
\end{prop}
In the proposition above, the union is that of subobjects of the codomain,
which is the same as the pushout over the intersection $\Lambda^{e}[S]\otimes\partial\Omega[T]$. 

\begin{proof}
(of Theorem 9.1 from Proposition 9.2) The theorem states that for
any tree $S$ and any inner edge $e\in S$, any map of dendroidal
sets\[
\varphi:\Lambda^{e}[S]\otimes X\rightarrow K\]
extends to some map (uniquely in the strict case)\[
\psi:\Omega[S]\otimes X\rightarrow K.\]
 By writing $X$ as the union of its skeleta, \[
X=\varinjlim Sk_{n}(X)\]
as in Section 4, and using the fact that the skeletal filtration is normal, we can build this extension $\psi$ by induction
on $n$. For $n=0$, $Sk_{0}(X)$ is a sum of copies of $\Omega[\eta]$,
the unit for the tensor product, so obviously the restriction $\varphi_{0}:\Lambda^{e}[S]\otimes Sk_{0}(X)\rightarrow K$
extends to a map \[
\psi_{0}:\Omega[S]\otimes Sk_{0}(X)\rightarrow K.\]
Suppose now that we have found an extension $\psi_{n}:\Omega[S]\otimes Sk_{n}(X)\rightarrow K$
of the restriction $\varphi_{n}:\Lambda^{e}[S]\otimes Sk_{n}(X)\rightarrow K$.
Consider the following diagram:\[
\xyC{5pt}\xymatrix{*++{\coprod\Lambda^{e}[S]\otimes\partial\Omega[T]}\ar[rr]\ar@{>->}[dd]\ar[rd] &  & *++{\coprod\Lambda^{e}[S]\otimes\Omega[T]}\ar@{>->}'[d][dd]\ar[rd]\\
 & *++{\Lambda^{e}[S]\otimes Sk_{n}(X)}\ar[rr]\ar@{>->}[dd] &  & *++{\Lambda^{e}[S]\otimes Sk_{n+1}(X)}\ar@{>->}[dd]\\
\coprod\Omega[S]\otimes\partial\Omega[T]\ar[rd]\ar'[r][rr] &  & \coprod\Omega[S]\otimes\Omega[T]\ar[rd]\\
 & \Omega[S]\otimes Sk_{n}(X)\ar[rr] &  & \Omega[S]\otimes Sk_{n+1}(X)}
\]
In this diagram, the top and bottom faces are pushouts given by the normal skeletal filtration of $X$. Now inscribe
the pushouts $U$ and $V$ in the back and front face, fitting into
a square\[
\xyC{15pt}\xymatrix{*+++{U}\ar@{>->}[d]\ar[r] & *+++{\coprod\Omega[S]\otimes\Omega[T]}\ar@{>->}[d]\\
V\ar[r] & \Omega[S]\otimes Sk_{n+1}(X)}
\]
The maps $\psi_{n}:\Omega[S]\otimes Sk_{n}(X)\rightarrow K$ and $\varphi_{n+1}:\Lambda^{e}[S]\otimes Sk_{n+1}(X)\rightarrow K$
together define a map $V\rightarrow K$. So, to find $\psi_{n+1},$
it suffices to prove that \[
\xymatrix{*++{V}\ar@{>->}[r] & \Omega[S]\otimes Sk_{n+1}(X)}
\]
 is anodyne. But, by a diagram chase argument, the square above is
a pushout, so in fact, it suffices to prove that $\xymatrix{*++{U}\ar@{>->}[r] & \coprod\Omega[S]\otimes\Omega[T]}
$ is anodyne. The latter map is a sum of copies of anodyne extensions
as in the statement of the proposition. 
\end{proof}
\begin{cor}
The monoidal structure on the category of coloured operads given by
the Boardman-Vogt tensor product is closed \emph{(}see \cite{DenSet}\emph{)}. It
is related to the closed monoidal structure on dendroidal sets by
two natural isomorphisms\[
\tau_{d}(N_{d}\mathcal{P}\otimes N_{d}\mathcal{Q})=\mathcal{P}\otimes_{BV}\mathcal{Q}\]
and \[
N_{d}(\underline{Hom}(\mathcal{Q},\mathcal{R}))=\underline{Hom}(N_{d}\mathcal{Q},N_{d}\mathcal{R})\]
 for any operads $\mathcal{P},\mathcal{Q}$ and $\mathcal{R}$. 
\end{cor}
\begin{proof}
The first isomorphism was proved in \cite{DenSet}. The second isomorphism
follows from the first one together with (the strict version of) Theorem
9.1, Theorem 6.1, and the fact that $N_{d}$ is fully faithful. 
\end{proof}
In the rest of this section, we will be concerned with the proof of
Proposition 9.2, and we fix $S,T$ and $e$ as in the statement of
the proposition from now on. Our strategy will be as follows. First,
let us write \[
A_{0}\subseteq\Omega[S]\otimes\Omega[T]\]
for the dendroidal set given by the image of $\Lambda^{e}[S]\otimes\Omega[T]\cup\Omega[S]\otimes\partial\Omega[T]$.
We are going to construct a sequence of dendroidal subsets \[
A_{0}\subseteq A_{1}\subseteq A_{2}\subseteq\cdots\subseteq A_{N}=\Omega[S]\otimes\Omega[T]\]
such that each inclusion is an anodyne extension. This will be done
by writing $\Omega[S]\otimes\Omega[T]$ as a union of representables,
as follows. We will explicitly describe a sequence of trees\[
T_{1},T_{2},\cdots,T_{N}\]
together with canonical monomorphisms (all called)\[
\xymatrix{*++{m:\Omega[T_{i}]}\ar@{>->}[r] & \Omega[S]\otimes\Omega[T]}
,\]
and we will write $m(T_{i})\subseteq\Omega[S]\otimes\Omega[T]$ for
the dendroidal subset given by the image of this monomorphism. We
will then define\[
\begin{array}{cccc}
 & A_{i+1}=A_{i}\cup m(T_{i+1}) & \quad\quad\quad\quad & (i=0,\cdots,N-1)\end{array}\]
and prove that each $\xymatrix{*++{A_{i}}\ar@{>->}[r] & A_{i+1}}
$ thus constructed is anodyne. For the rest of this section, we will
fix planar structures on the trees $S$ and $T$. These will then
induce a natural planar structure on each of the trees $T_{i}$, and
avoid unnecessary discussion involving automorphisms in the category
$\Omega$. 

To define the $T_{i}$, let us think of the vertices of $S$ as \emph{white}
(drawn $\circ$) and those of $T$ as \emph{black} (drawn $\bullet$).
The edges of $T_{i}$ are (labelled by) pairs $(a,x)$ where $a$
is an edge of $S$ and $x$ one of $T$. We refer to $a$ as the $S$-\emph{colour}
of this edge $(a,x)$, and to $x$ as its $T$-\emph{colour}. There
are two kinds of vertices in $T_{i}$ (corresponding to the generators
for $\Omega[S]\otimes\Omega[T]$ coming from vertices of $S$ or of
$T$). There are \emph{white} vertices in $T_{i}$ labelled 
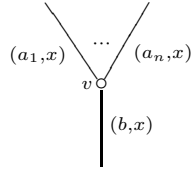
\[
\xymatrix{*{}\ar@{-}[dr]_{(a_{1},x)} &  & *{}\ar@{-}[dl]^{(a_{n},x)}\\
\ar@{}[r]|{\quad v} & *{\circ}\ar@{}[u]|{\cdots}\ar@{-}[d]^{(b,x)}\\
 & *{}}
\]
where $v$ is a vertex in $S$ with input edges $a_{1},\cdots,a_{n}$
and output edge $b$, while $x$ is an edge of $T$; and there are
\emph{black} vertices in $T_{i}$ labelled
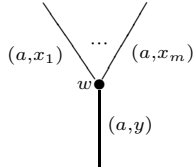
\[
\xymatrix{*{}\ar@{-}[dr]_{(a,x_{1})} &  & *{}\ar@{-}[dl]^{(a,x_{m})}\\
\ar@{}[r]|{\quad w} & *{\bullet}\ar@{}[u]|{\cdots}\ar@{-}[d]^{(a,y)}\\
 & *{}}
\]
where $w$ is a vertex in $T$ with input edges $x_{1},\cdots,x_{m}$
and output edge $y$, while $a$ is an edge in $S$. Moreover, each
such tree $T_{i}$ is \emph{maximal,} in the sense that its output
(root) edge is labelled $(r_{S},r_{T})$ where $r_{S}$ and $r_{T}$
are the roots of $S$ and $T$, and its input edges are labelled by
all pairs $(a,x)$ where $a$ is an input edge of $S$ and $x$ one
of $T$. 

All the possible such trees $T_{i}$ come in a natural (partial) order.
The minimal tree $T_{1}$ in the poset is the one obtained by stacking
a copy of the black tree $T$ on top of each of the input edges of
the white tree $S$. Or, more precisely, on the bottom of $T_{1}$
there is a copy $S\otimes r_{T}$ of the tree $S$ all whose edges
are renamed $(a,r_{T})$ where $r_{T}$ is the output edge at the
root of $T$. For each input edge $b$ of $S$, a copy of $T$ is
grafted on the edge $(b,r)$ of $S\otimes r$, with edges $x$ in
$T$ renamed $(b,x)$. The maximal tree $T_{N}$ in the poset is the
similar tree with copies of the white tree $S$ grafted on each of
the input edges of the black tree. Pictorially $T_{1}$ looks like 

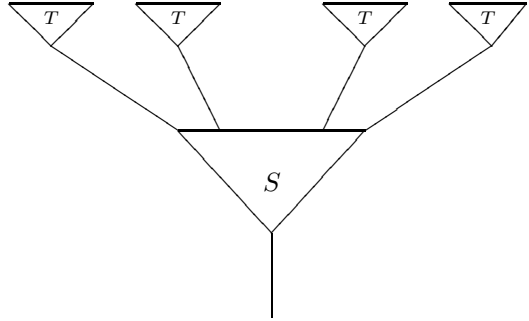
\[
\xyR{10pt}\xyC{10pt}\xymatrix{*{}\ar@{-}[rr]_{T}\ar@{-}[rd] &  & *{}\ar@{-}[dl] & *{}\ar@{-}[rr]_{T}\ar@{-}[rd] &  & *{}\ar@{-}[dl] &  & *{}\ar@{-}[rr]_{T}\ar@{-}[rd] &  & *{}\ar@{-}[dl] & *{}\ar@{-}[rr]_{T}\ar@{-}[rd] &  & *{}\ar@{-}[dl]\\
 & *{}\ar@{-}[ddrrr] &  &  & *{}\ar@{-}[ddr] &  &  &  & *{}\ar@{-}[ddl] &  &  & *{}\ar@{-}[ddlll]\\
\\ &  &  &  & *{}\ar@{-}[rrrr]\ar@{-}[ddrr] & *{} &  & *{} & *{}\ar@{-}[ddll]\\
 &  &  &  &  &  & S\\
 &  &  &  &  &  & *{}\ar@{-}[dd]\\
\\ &  &  &  &  &  & *{}}
\]
and $T_{N}$ looks like 

\[
\xyR{10pt}\xyC{10pt}
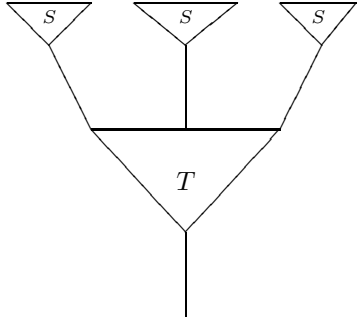
\xymatrix{*{}\ar@{-}[dr]\ar@{-}[rr]_{S} &  & *{}\ar@{-}[dl] & *{}\ar@{-}[dr]\ar@{-}[rr]_{S} &  & *{}\ar@{-}[dl] & *{}\ar@{-}[dr]\ar@{-}[rr]_{S} &  & *{}\ar@{-}[dl]\\
 & *{}\ar@{-}[ddr] &  &  & *{}\ar@{-}[dd] &  &  & *{}\ar@{-}[ddl]\\
\\ &  & *{}\ar@{-}[rrrr]\ar@{-}[rrdd] &  & *{} &  & *{}\ar@{-}[ddll]\\
 &  &  &  & T\\
 &  &  &  & *{}\ar@{-}[dd]\\
\\ &  &  &  & *{}}
\]

The intermediate trees $T_{k}$ ($1<k<N)$ are obtained by letting
the black vertices in $T_{1}$ slowly percolate in all possible ways
towards the root of the tree. Each $T_{k}$ is obtained from an earlier
$T_{l}$ by replacing a configuration

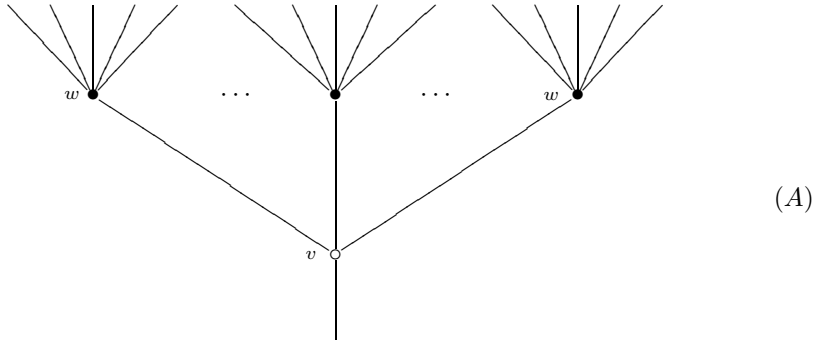
\[
\xyR{10pt}\xyC{10pt}\xymatrix{ & *{}\ar@{-}[ddrr] & *{}\ar@{-}[ddr] & *{}\ar@{-}[dd] & *{}\ar@{-}[ddl] & *{}\ar@{-}[ddll] & *{}\ar@{-}[ddrr] & *{}\ar@{-}[ddr] & *{}\ar@{-}[dd] & *{}\ar@{-}[ddl] & *{}\ar@{-}[ddll] & *{}\ar@{-}[ddrr] & *{}\ar@{-}[ddr] & *{}\ar@{-}[dd] & *{}\ar@{-}[ddl] & *{}\ar@{-}[ddll] & *{} & \,\\
\\ &  & *{}\ar@{}[r]|{w} & *{\bullet}\ar@{-}[rrrrrddd] &  &  & \cdots &  & *{\bullet}\ar@{-}[ddd] &  & \cdots &  & *{}\ar@{}[r]|{w\,\,\,} & *{\bullet}\ar@{-}[dddlllll]\\
\\\quad &  &  &  &  &  &  &  &  &  &  &  &  &  &  &  &  & \quad(A)\\
 &  &  &  &  &  &  & \ar@{}[r]|{v\,\,} & *{\circ}\ar@{-}[dd]\\
\\ &  &  &  &  &  &  &  & *{}}
\]
 in $T_{l}$ by \[
\xyR{10pt}\xyC{10pt}
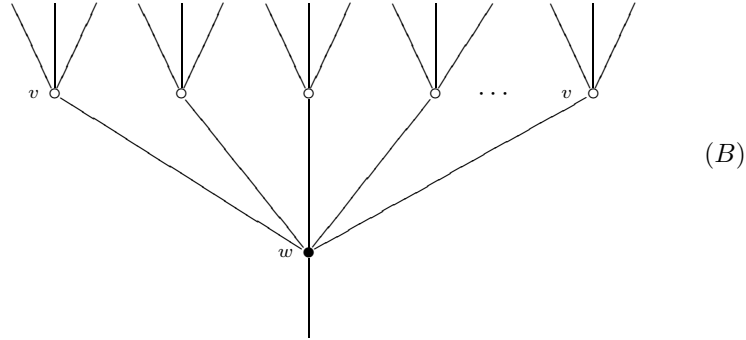
\xymatrix{ & *{}\ar@{-}[ddr] & *{}\ar@{-}[dd] & *{}\ar@{-}[ddl] & *{}\ar@{-}[ddr] & *{}\ar@{-}[dd] & *{}\ar@{-}[ddl] & *{}\ar@{-}[ddr] & *{}\ar@{-}[dd] & *{}\ar@{-}[ddl] & *{}\ar@{-}[ddr] & *{}\ar@{-}[dd] & *{}\ar@{-}[ddl] & *{}\ar@{-}[ddr] & *{}\ar@{-}[dd] & *{}\ar@{-}[ddl] & \,\\
\\ & *{}\ar@{}[r]|{v} & *{\circ}\ar@{-}[dddrrrrrr] &  &  & *{\circ}\ar@{-}[dddrrr] &  &  & *{\circ}\ar@{-}[ddd] &  &  & *{\circ}\ar@{-}[dddlll] & \cdots & *{}\ar@{}[r]|{v\,\,\,} & *{\circ}\ar@{-}[dddllllll]\\
\quad &  &  &  &  &  &  &  &  &  &  &  &  &  &  &  & \quad(B)\\
\\ &  &  &  &  &  &  & *{}\ar@{}[r]|{w\,} & *{\bullet}\ar@{-}[dd]\\
\\ &  &  &  &  &  &  &  & *{}}
\]
in $T_{k}$. More explicitly, if $v$ and $w$ are vertices in $S$
and $T$, \[
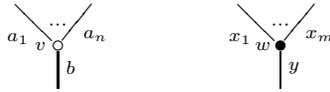
\begin{array}{ccc}
\xyC{10pt}\xyR{10pt}\xymatrix{*{}\ar@{-}[dr]_{a_{1}} &  & *{}\ar@{-}[dl]^{a_{n}}\\
\ar@{}[r]|{\,\, v} & *{\circ}\ar@{}[u]|{\cdots}\ar@{-}[d]^{b}\\
 & *{}}
 & \quad\quad & \xymatrix{*{}\ar@{-}[dr]_{x_{1}} &  & *{}\ar@{-}[dl]^{x_{m}}\\
\ar@{}[r]|{\,\, w} & *{\bullet}\ar@{}[u]|{\cdots}\ar@{-}[d]^{y}\\
 & *{}}
\end{array}\]
then the edges in $(A)$ are named 
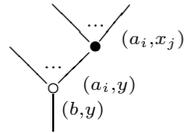
\[
\xymatrix{ & *{}\ar@{-}[dr] &  & *{}\ar@{-}[dl]^{(a_{i},x_{j})}\\
*{}\ar@{-}[dr] &  & *{\bullet}\ar@{-}[dl]^{(a_{i},y)}\ar@{}[u]|{\cdots}\\
 & *{\circ}\ar@{-}[d]^{(b,y)}\ar@{}[u]|{\cdots}\\
 & *{}}
\]
and those in $(B)$ are named
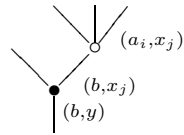
\[
\xymatrix{ & *{}\ar@{-}[dr] & *{}\ar@{-}[d] & *{}\ar@{-}[dl]^{(a_{i},x_{j})}\\
*{}\ar@{-}[dr] &  & *{\circ}\ar@{-}[dl]^{(b,x_{j})}\\
 & *{\bullet}\ar@{-}[d]^{(b,y)}\\
 & *{}}
\]
We will refer to these trees $T_{i}$ as the \emph{percolation schemes}
for $S$ and $T$, and if $T_{k}$ is obtained from $T_{l}$ by replacing
$(A)$ by $(B)$, then we will say that $T_{l}$ is obtained by a
\emph{single percolation step}.

\begin{example}
Many of the typical phenomena that we will encounter already occur
for the following two trees $S$ and $T$; here, we have singled out
one particular edge $e$ in $S$, we've numbered the edges of $T$
as $1,\cdots,5$, and denoted the colour $(e,i)$ in $T_{i}$ by $e_{i}$.
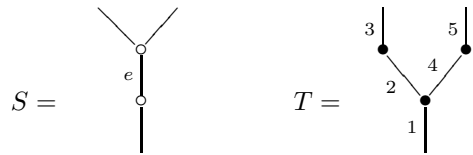
\[
\begin{array}{ccc}
\xymatrix{ & *{\,}\ar@{-}[dr] &  & *{\,}\ar@{-}[dl]\\
 &  & *{\circ}\ar@{-}[d]_{e}\\
S= &  & *{\circ}\ar@{-}[d]\\
 &  & *{\,}}
 & \quad\quad & \xymatrix{ & *{\,}\ar@{-}[d]_{3} &  & *{\,}\ar@{-}[d]_{5}\\
 & *{\bullet}\ar@{-}[dr]_{2} &  & *{\bullet}\ar@{-}[dl]_{4}\\
T= &  & *{\bullet}\ar@{-}[d]_{1}\\
 &  & *{\,}}
\end{array}\]

There are 14 percolation schemes $T_{1},\cdots,T_{14}$ in this case.
Here is the complete list of them:

$\xyR{5pt}\xyC{5pt}
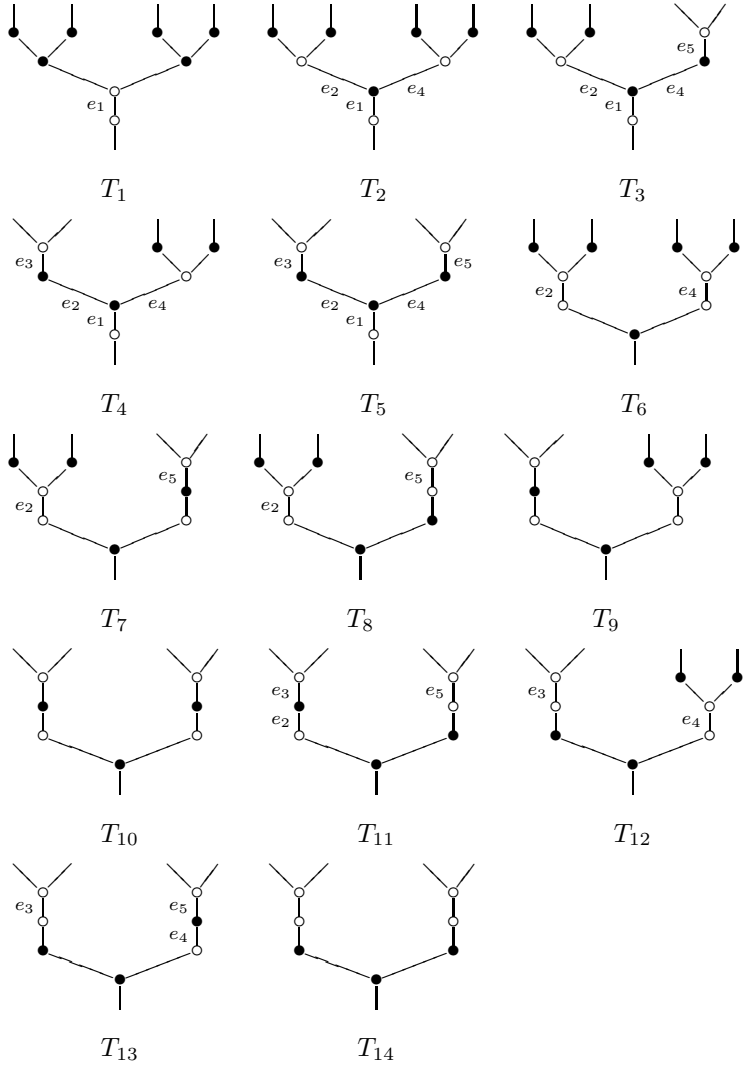
\xymatrix{*{}\ar@{-}[d] &  & *{}\ar@{-}[d] &  & *{}\ar@{-}[d] &  & *{}\ar@{-}[d]\\
*{\bullet}\ar@{-}[dr] &  & *{\bullet}\ar@{-}[dl] &  & *{\bullet}\ar@{-}[dr] &  & *{\bullet}\ar@{-}[dl]\\
 & *{\bullet}\ar@{-}[drr] &  &  &  & *{\bullet}\ar@{-}[dll]\\
 &  &  & *{\circ}\ar@{-}[d]_{e_{1}}\\
 &  &  & *{\circ}\ar@{-}[d]\\
 &  &  & *{}\\
 &  &  & T_{1}}
$~~~~~$\xyR{5pt}\xyC{5pt}\xymatrix{*{}\ar@{-}[d] &  & *{}\ar@{-}[d] &  & *{}\ar@{-}[d] &  & *{}\ar@{-}[d]\\
*{\bullet}\ar@{-}[dr] &  & *{\bullet}\ar@{-}[dl] &  & *{\bullet}\ar@{-}[dr] &  & *{\bullet}\ar@{-}[dl]\\
 & *{\circ}\ar@{-}[drr]_{e_{2}} &  &  &  & *{\circ}\ar@{-}[dll]^{e_{4}}\\
 &  &  & *{\bullet}\ar@{-}[d]_{e_{1}}\\
 &  &  & *{\circ}\ar@{-}[d]\\
 &  &  & *{}\ar@{-}[]\\
 &  &  & T_{2}}
$~~~~~$\xyR{5pt}\xyC{5pt}\xymatrix{*{}\ar@{-}[d] &  & *{}\ar@{-}[d] &  & *{}\ar@{-}[dr] &  & *{}\ar@{-}[dl]\\
*{\bullet}\ar@{-}[dr] &  & *{\bullet}\ar@{-}[dl] &  &  & *{\circ}\ar@{-}[d]_{e_{5}}\\
 & *{\circ}\ar@{-}[drr]_{e_{2}} &  &  &  & *{\bullet}\ar@{-}[dll]^{e_{4}}\\
 &  &  & *{\bullet}\ar@{-}[d]_{e_{1}}\\
 &  &  & *{\circ}\ar@{-}[d]\\
 &  &  & *{}\\
 &  &  & T_{3}}
$

$\xyR{5pt}\xyC{5pt}\xymatrix{*{}\ar@{-}[dr] &  & *{}\ar@{-}[dl] &  & *{}\ar@{-}[d] &  & *{}\ar@{-}[d]\\
 & *{\circ}\ar@{-}[d]_{e_{3}} &  &  & *{\bullet}\ar@{-}[dr] &  & *{\bullet}\ar@{-}[dl]\\
 & *{\bullet}\ar@{-}[drr]_{e_{2}} &  &  &  & *{\circ}\ar@{-}[dll]^{e_{4}}\\
 &  &  & *{\bullet}\ar@{-}[d]_{e_{1}}\\
 &  &  & *{\circ}\ar@{-}[d]\\
 &  &  & *{}\\
 &  &  & T_{4}}
$~~~~~$\xyR{5pt}\xyC{5pt}\xymatrix{*{}\ar@{-}[dr] &  & *{}\ar@{-}[dl] &  & *{}\ar@{-}[dr] &  & *{}\ar@{-}[dl]\\
 & *{\circ}\ar@{-}[d]_{e_{3}} &  &  &  & *{\circ}\ar@{-}[d]^{e_{5}}\\
 & *{\bullet}\ar@{-}[drr]_{e_{2}} &  &  &  & *{\bullet}\ar@{-}[dll]^{e_{4}}\\
 &  &  & *{\bullet}\ar@{-}[d]_{e_{1}}\\
 &  &  & *{\circ}\ar@{-}[d]\\
 &  &  & *{}\\
 &  &  & T_{5}}
$~~~~~$\xyR{5pt}\xyC{5pt}\xymatrix{*{}\ar@{-}[d] &  & *{}\ar@{-}[d] &  & *{}\ar@{-}[d] &  & *{}\ar@{-}[d]\\
*{\bullet}\ar@{-}[dr] &  & *{\bullet}\ar@{-}[dl] &  & *{\bullet}\ar@{-}[dr] &  & *{\bullet}\ar@{-}[dl]\\
 & *{\circ}\ar@{-}[d]_{e_{2}} &  &  &  & *{\circ}\ar@{-}[d]_{e_{4}}\\
 & *{\circ}\ar@{-}[drr] &  &  &  & *{\circ}\ar@{-}[dll]\\
 &  &  & *{\bullet}\ar@{-}[d]\\
 &  &  & *{}\\
 &  &  & T_{6}}
$

$\xyR{5pt}\xyC{5pt}\xymatrix{*{}\ar@{-}[d] &  & *{}\ar@{-}[d] &  & *{}\ar@{-}[dr] &  & *{}\ar@{-}[dl]\\
*{\bullet}\ar@{-}[dr] &  & *{\bullet}\ar@{-}[dl] &  &  & *{\circ}\ar@{-}[d]_{e_{5}}\\
 & *{\circ}\ar@{-}[d]_{e_{2}} &  &  &  & *{\bullet}\ar@{-}[d]\\
 & *{\circ}\ar@{-}[drr] &  &  &  & *{\circ}\ar@{-}[dll]\\
 &  &  & *{\bullet}\ar@{-}[d]\\
 &  &  & *{}\\
 &  &  & T_{7}}
$~~~~~$\xyR{5pt}\xyC{5pt}\xymatrix{*{}\ar@{-}[d] &  & *{}\ar@{-}[d] &  & *{}\ar@{-}[dr] &  & *{}\ar@{-}[dl]\\
*{\bullet}\ar@{-}[dr] &  & *{\bullet}\ar@{-}[dl] &  &  & *{\circ}\ar@{-}[d]_{e_{5}}\\
 & *{\circ}\ar@{-}[d]_{e_{2}} &  &  &  & *{\circ}\ar@{-}[d]\\
 & *{\circ}\ar@{-}[drr] &  &  &  & *{\bullet}\ar@{-}[dll]\\
 &  &  & *{\bullet}\ar@{-}[d]\\
 &  &  & *{}\\
 &  &  & T_{8}}
$~~~~~$\xyC{5pt}\xyR{5pt}\xymatrix{*{}\ar@{-}[dr] &  & *{}\ar@{-}[dl] &  & *{}\ar@{-}[d] &  & *{}\ar@{-}[d]\\
 & *{\circ}\ar@{-}[d] &  &  & *{\bullet}\ar@{-}[dr] &  & *{\bullet}\ar@{-}[dl]\\
 & *{\bullet}\ar@{-}[d] &  &  &  & *{\circ}\ar@{-}[d]\\
 & *{\circ}\ar@{-}[drr] &  &  &  & *{\circ}\ar@{-}[dll]\\
 &  &  & *{\bullet}\ar@{-}[d]\\
 &  &  & *{}\\
 &  &  & T_{9}}
$

$\xyC{5pt}\xyR{5pt}\xymatrix{*{}\ar@{-}[dr] &  & *{}\ar@{-}[dl] &  & *{}\ar@{-}[dr] &  & *{}\ar@{-}[dl]\\
 & *{\circ}\ar@{-}[d] &  &  &  & *{\circ}\ar@{-}[d]\\
 & *{\bullet}\ar@{-}[d] &  &  &  & *{\bullet}\ar@{-}[d]\\
 & *{\circ}\ar@{-}[drr] &  &  &  & *{\circ}\ar@{-}[dll]\\
 &  &  & *{\bullet}\ar@{-}[d]\\
 &  &  & *{}\\
 &  &  & T_{10}}
$~~~~~$\xyR{5pt}\xyC{5pt}\xymatrix{*{}\ar@{-}[dr] &  & *{}\ar@{-}[dl] &  & *{}\ar@{-}[dr] &  & *{}\ar@{-}[dl]\\
 & *{\circ}\ar@{-}[d]_{e_{3}} &  &  &  & *{\circ}\ar@{-}[d]_{e_{5}}\\
 & *{\bullet}\ar@{-}[d]_{e_{2}} &  &  &  & *{\circ}\ar@{-}[d]\\
 & *{\circ}\ar@{-}[drr] &  &  &  & *{\bullet}\ar@{-}[dll]\\
 &  &  & *{\bullet}\ar@{-}[d]\\
 &  &  & *{}\\
 &  &  & T_{11}}
$~~~~~$\xyC{5pt}\xyR{5pt}\xymatrix{*{}\ar@{-}[dr] &  & *{}\ar@{-}[dl] &  & *{}\ar@{-}[d] &  & *{}\ar@{-}[d]\\
 & *{\circ}\ar@{-}[d]_{e_{3}} &  &  & *{\bullet}\ar@{-}[dr] &  & *{\bullet}\ar@{-}[dl]\\
 & *{\circ}\ar@{-}[d] &  &  &  & *{\circ}\ar@{-}[d]_{e_{4}}\\
 & *{\bullet}\ar@{-}[drr] &  &  &  & *{\circ}\ar@{-}[dll]\\
 &  &  & *{\bullet}\ar@{-}[d]\\
 &  &  & *{}\\
 &  &  & T_{12}}
$

$\xyC{5pt}\xyR{5pt}\xymatrix{*{}\ar@{-}[dr] &  & *{}\ar@{-}[dl] &  & *{}\ar@{-}[dr] &  & *{}\ar@{-}[dl]\\
 & *{\circ}\ar@{-}[d]_{e_{3}} &  &  &  & *{\circ}\ar@{-}[d]_{e_{5}}\\
 & *{\circ}\ar@{-}[d] &  &  &  & *{\bullet}\ar@{-}[d]_{e_{4}}\\
 & *{\bullet}\ar@{-}[drr] &  &  &  & *{\circ}\ar@{-}[dll]\\
 &  &  & *{\bullet}\ar@{-}[d]\\
 &  &  & *{}\\
 &  &  & T_{13}}
$~~~~~$\xyC{5pt}\xyR{5pt}\xymatrix{*{}\ar@{-}[dr] &  & *{}\ar@{-}[dl] &  & *{}\ar@{-}[dr] &  & *{}\ar@{-}[dl]\\
 & *{\circ}\ar@{-}[d] &  &  &  & *{\circ}\ar@{-}[d]\\
 & *{\circ}\ar@{-}[d] &  &  &  & *{\circ}\ar@{-}[d]\\
 & *{\bullet}\ar@{-}[drr] &  &  &  & *{\bullet}\ar@{-}[dll]\\
 &  &  & *{\bullet}\ar@{-}[d]\\
 &  &  & *{}\\
 &  &  & T_{14}}
$
\end{example}
As claimed, there is a partial order on the percolation schemes $T_{1},\cdots,T_{N}$
for $S\otimes T$, in which $T_{1}$ (copies of $T$ on top of $S$)
is the minimal element and $T_{N}$ (copies of $S$ on top of $T$)
the maximal one. The partial order is given by defining $T\le T'$
whenever the percolation scheme $T'$ can be obtained from the percolation
scheme $T$ by a sequence of percolations. For example, the poset
structure on the percolation trees above is:\[
\xyC{10pt}\xyR{10pt}
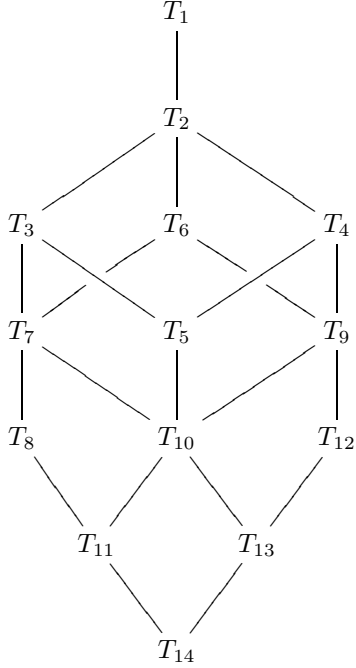
\xymatrix{ &  & T_{1}\ar@{-}[dd]\\
\\ &  & T_{2}\ar@{-}[ddrr]\ar@{-}[dd]\ar@{-}[ddll]\\
\\T_{3}\ar@{-}[dd]\ar@{-}[ddrr] &  & T_{6}\ar@{-}'[dl][ddll]\ar@{-}'[dr][ddrr] &  & T_{4}\ar@{-}[dd]\ar@{-}[ddll]\\
 & \, &  & \,\\
T_{7}\ar@{-}[ddrr]\ar@{-}[dd] &  & T_{5}\ar@{-}[dd] &  & T_{9}\ar@{-}[ddll]\ar@{-}[dd]\\
 & \,\\
T_{8}\ar@{-}[ddr] &  & T_{10}\ar@{-}[ddr]\ar@{-}[ddl] &  & T_{12}\ar@{-}[ddl]\\
\\ & T_{11}\ar@{-}[ddr] &  & T_{13}\ar@{-}[ddl]\\
\\ &  & T_{14}}
\]
The planar structures of $S$ and $T$ provide a way to refine this
partial order by a linear order. It is not important exactly how this
is done, but we shall from now on assume that the percolation schemes
for $S$ and $T$ are ordered $T_{1},\cdots,T_{N}$ where $T_{i}$
comes before $T_{j}$ only if $T_{i}\le T_{j}$ in the partial order. 

\begin{lem}
(and notation) Each percolation scheme $T_{i}$ is equipped with a
canonical monomorphism \[
\xymatrix{*++{m:\Omega[T_{i}]}\ar@{>->}[r] & \Omega[S]\otimes\Omega[T]}
.\]
 The dendroidal subset given by the image of this monomorphism will
be denoted \[
m(T_{i})\subseteq\Omega[S]\otimes\Omega[T].\]
 
\end{lem}
\begin{proof}
The vertices of the dendroidal set $\Omega[T_{i}]$ are the edges
of the tree $T_{i}$. The map $m$ is completely determined by asking
it to map an edge named $(a,x)$ in $T_{i}$ to the vertex with the
same name in $\Omega[S]\otimes\Omega[T]$. This map is a monomorphism.
In fact, any map \[
\Omega[R]\rightarrow X,\]
 from a representable dendroidal set to an arbitrary one, is a monomorphism
as soon as the map $\Omega[R]_{\eta}\rightarrow X_{\eta}$ on vertices
is. 
\end{proof}
Before we continue, we need to introduce a bit of terminology for
trees, i.e., for objects of $\Omega$. Let $R$ be such a tree. A
map $R'\rightarrow R$ which is a composition of basic face maps (maps
of type $(ii)$ or $(iii)$ in Section 3) will also be referred to
as a \emph{face} of $R$, just like for simplicial sets. If it is
a composition of \emph{inner} faces (resp. \emph{outer} faces), the
map $\xymatrix{*++{R'}\ar@{>->}[r] & R}
$ will be called an \emph{inner face} (resp. \emph{outer face)} of
$R$. A \emph{top face} of $R$ is an outer face map $\partial_{v}:R'\rightarrow R$
where $R'$ is obtained by deleting a top vertex from $R$. An \emph{initial
segment} $\xymatrix{*++{R'}\ar@{>->}[r] & R}
$ is a composition of top faces (it is a special kind of outer face
of $R$). If $v$ is the vertex above the root of $R$ and $e$ is
an input edge of $v$, then $R$ contains a subtree $R'$ whose root
is $e$. We'll refer to an inclusion of this kind as a \emph{bottom} face
of $R$ (it is again a special kind of outer face). In all these cases,
we'll often leave the monomorphism $\xymatrix{*++{R'}\ar@{>->}[r] & R}
$ implicit, and apply the same terminology not only to the map $\xymatrix{*++{R'}\ar@{>->}[r] & R}
$ but also to the tree $R'$. 

For example, for the tree $T$ constructed above 
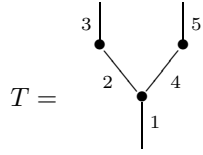
\[
\xymatrix{ & *{\,}\ar@{-}[d]_{3} &  & *{\,}\ar@{-}[d]^{5}\\
 & *{\bullet}\ar@{-}[dr]_{2} &  & *{\bullet}\ar@{-}[dl]^{4}\\
T= &  & *{\bullet}\ar@{-}[d]^{1}\\
 &  & *{\,}}
\]
 The following sub-trees are examples of, respectively, a top face,
an initial segment, a bottom face, and an inner face:\[
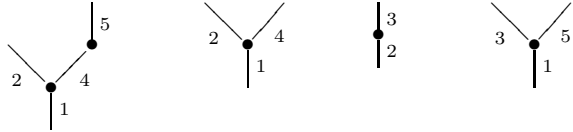
\begin{array}{ccccccc}
\xyC{10pt}\xyR{10pt}\xymatrix{ &  & *{\,}\ar@{-}[d]^{5}\\
*{\,}\ar@{-}[dr]_{2} &  & *{\bullet}\ar@{-}[dl]^{4}\\
 & *{\bullet}\ar@{-}[d]^{1}\\
 & *{\,}}
 & \quad & \xymatrix{*{\,}\ar@{-}[dr]_{2} &  & *{\,}\ar@{-}[dl]^{4}\\
 & *{\bullet}\ar@{-}[d]^{1}\\
 & *{\,}}
 & \quad & \xymatrix{*{\,}\ar@{-}[d]^{3}\\
*{\bullet}\ar@{-}[d]^{2}\\
*{\,}}
 & \quad & \xymatrix{*{\,}\ar@{-}[dr]_{3} &  & *{\,}\ar@{-}[dl]^{5}\\
 & *{\bullet}\ar@{-}[d]^{1}\\
 & *{\,}}
\end{array}\]

\begin{rem}
We observe the following simple properties, which we will repeatedly
use in the proofs of the lemmas below. In stating these properties
and below, we denote by $m(R)$ the image of the composition of the
inclusion $\xymatrix{*++{\Omega[R]}\ar@{>->}[r] & \Omega[T_{i}]}
$ given by a subtree (a face) $R$ of $T_{i}$ and the canonical monomorphism
\[
\xyC{15pt}\xyR{15pt}\xymatrix{*++{m:\Omega[T_{i}]}\ar@{>->}[r] & \Omega[S]\otimes\Omega[T]}
.\]

(i) Let $R$ be a subtree of $T_{i}$. If $m(R)\subseteq A_{0}$ then
$R$ misses a $T$-colour, or an $S$-colour other than $e$, or a
stump of either $S$ or $T$. Here, a stump is a top vertex of valence
zero (i.e., without input edges). We say that $R$ ''misses'' such
a stump $v\in S$, for example, if $m(R)\subseteq\partial_{v}[S]\otimes\Omega[T]$.
The tree $R$ is a sub-tree of $T_{i}$, where edges
are coloured by pairs $(a,x)$, where $a$ is an $S$-colour and $x$
a $T$-colour. By saying that $R$ ''misses'' a $T$-colour $y$,
we mean that none of the colours $(a,x)$ occurring in $R$ has $x=y$
as second coordinate. ''Missing an $S$-colour'' is interpreted
similarly. 

(ii) This implies in particular that for any bottom face $\xymatrix{R\ar[r] & T_{i}}
$ of any percolation scheme $T_{i}$ the dendroidal set $m(R)$ is
contained in $A_{0}$, because it must miss either the root colour
$r_{S}$ (in case the root of $T_{i}$ is white), or the root colour
$r_{T}$ (in case the root of $T_{i}$ is black), and $r_{S}\ne e$
because $e$ is assumed inner.

(iii) If $F,G$ are faces of $T_{i}$ then $F$ is a face of $G$
iff $m(F)\subseteq m(G)$. (This is clear from the fact that the map
from $\Omega[R]$ onto its image $m(R)$ is an isomorphism of dendroidal
sets.) 

(iv) Let $Q$ and $R$ be initial segments of $T_{i}$, and let $F$
be an inner face of $Q$. If $m(F)\subseteq m(R)$ then also $m(Q)\subseteq m(R)$
(and hence $Q$ is a face of $R$, by (iii)). In fact, let $Inn(Q)$
denote the set of all inner edges of $Q$ and $\xymatrix{*++{Q/Inn(Q)}\ar@{>->}[r] & Q}
$ the inner face of $Q$ given by contracting all these. Then if $m(Q/Inn(Q))\subseteq m(R/Inn(R))$,
it follows by comparing labels of input edges of $Q$ and $R$ that
$Q$ is a face of $R$. 
\end{rem}
These remarks prepare the ground for the following lemma. Recall that
$A_{k}=A_{0}\cup m(T_{1})\cup\cdots\cup m(T_{k})$, where $m(T_{i})$
is the image in $\Omega[S]\otimes\Omega[T]$ of the dendroidal set
$\Omega[T_{i}]$. 

\begin{lem}
Let $R,Q_{1},\cdots,Q_{p}$ be a family of initial segments in $T_{k+1}$
and write $B=m(Q_{1})\cup\cdots\cup m(Q_{p})\subseteq\Omega[S]\otimes\Omega[T]$.
Suppose 

(i) For every top face $F$ of $R$, $m(F)\subseteq A_{k}\cup B$. 

(ii) There exists an edge $\xi$ in $R$ such that for every inner
face $\xymatrix{*++{F}\ar@{>->}[r] & R}
$, if $m(F)$ is not contained in $A_{k}\cup B$ then neither is $m(F/(\xi))$.
\\
Then the inclusion $\xymatrix{*++{A_{k}\cup B}\ar@{>->}[r] & A_{k}\cup B\cup m(R)}
$ is anodyne.
\end{lem}
We call $\xi$ a characteristic edge of $R$ with respect to $Q_{1},\cdots,Q_{p}.$

\begin{proof}
If $m(R)\subseteq A_{k}\cup B$ there is nothing to prove. If not,
then by (ii), $m(R/(\xi))$ is not contained in $A_{k}\cup B$. Let
\[
\xi=\xi_{0},\xi_{1},\cdots,\xi_{n}\]
 be all the inner edges in $R$ such that the dendroidal set $m(R/(\xi_{i}))$
is not contained in $A_{k}\cup B$. For a sub-sequence $\xi_{i_{1}},\cdots,\xi_{i_{p}}$
of these $\xi_{0},\cdots,\xi_{n}$, we have the dendroidal subset
of $\Omega[S]\otimes\Omega[T]$,\[
\begin{array}{ccc}
\quad\quad\quad\quad\quad\quad\quad\quad\quad\quad\quad & m(R/(\xi_{i_{1}},\cdots,\xi_{i_{p}})), & \quad\quad\quad\quad\quad\quad\quad\quad\quad\quad\quad(1)\end{array}\]
 obtained by contracting each of $\xi_{i_{1}},\cdots,\xi_{i_{p}}$
and composing with $m:\Omega[T_{k+1}]\rightarrow\Omega[S]\otimes\Omega[T]$.
We are going to consider a sequence of anodyne extensions\[
\xymatrix{A_{k}\cup B=B_{0}\,\,\ar@{>->}[r] & B_{1}\,\,\ar@{>->}[r] & \cdots\,\,\ar@{>->}[r] & B_{2^{n}}=A_{k}\cup B\cup m(R)}
\]
 by considering images of faces of $\Omega[R]$ of this type (1).

Consider first \[
R_{(0)}=m(R/(\xi_{1},\cdots,\xi_{n})).\]
 If $m(R_{(0)})$ is contained in $A_{k}\cup B$, let $B_{1}=B_{0}=A_{k}\cup B$.
Otherwise, let $B_{1}$ be the pushout\[
\xyC{20pt}\xyR{20pt}\xymatrix{m(\Lambda^{\xi_{0}}R_{(0)})\ar[r]\ar[d] & B_{0}\ar[d]\\
m(R_{(0)})\ar[r] & B_{1}}
\]
 Notice that $m(\Lambda^{\xi_{0}}R_{(0)})$ is indeed contained in $B_{0}=A_{k}\cup B$.
For, any outer face $F$ of $R_{(0)}$ is a face of an outer face
$G$ of $R$\[
\xymatrix{F\ar[r]\ar[d] & R/(\xi_{1},\cdots,\xi_{n})=R_{(0)}\ar[d]\\
G\ar[r] & R}
\]
 if $G$ is a top face, then $m(G)\subseteq A_{k}\cup B$ by assumption
(i); and if $G$ is a bottom face, it already factors through $A_{0}\subseteq A_{k}$
(cf Remark 9.6 (ii) before the lemma). On the other hand, if $F\subseteq R_{(0)}$
is an inner face of $R_{(0)}$ given by contracting an edge $\zeta$
in $R/(\xi_{1},\cdots,\xi_{n})$, then $F$ is a face of $R/(\zeta)$.
So if $m(F)\nsubseteq B_{0}$ then $m(R/(\zeta))$ wouldn't be contained
in $B_{0}$ either, and hence $\zeta$ must be one of $\xi_{0},\cdots,\xi_{n}$.
But $\xi_{1},\cdots,\xi_{n}$ are no longer edges in $R/(\xi_{1},\cdots,\xi_{n})$,
so $\zeta$ must be $\xi_{0}$. This shows that for any inner face
$F$ of $R_{(0)}$ other then $R_{(0)}/(\xi_{0})$, the dendroidal
set $m(F)$ is contained in $B_{0}$, as claimed.

Next, consider all sub-sequences $(\xi_{1},\cdots,\hat{\xi}_{i},\cdots,\xi_{n})$,
and the faces \[
R_{(i)}=R/(\xi_{1},\cdots,\hat{\xi}_{i},\cdots\xi_{n})\,\,\,\,\,\, i=1,\cdots,n\]
 We will define \[
B_{2},\cdots,B_{n+1}\]
 by considering these $R_{(1)},\cdots,R_{(n)}$. Suppose $B_{1},\cdots,B_{i}$
have been defined. Consider $R_{(i)}$ to form $B_{i+1}$. If its
image $m(R_{(i)})$ is contained in $B_{i}$, let $B_{i+1}=B_{i}$.
Otherwise, $m(R_{(i)})\rightarrow\Omega[S]\otimes\Omega[T]$ does
not factor through $B_{i}$, and a fortiori doesn't factor through
$A_{k}\cup B=B_{0}$ either. So by assumption (ii), we have that $m(R_{(i)}/(\xi_{0}))\nsubseteq A_{k}\cup B$.
But then $m(R_{i}/(\xi_{0}))$ is not contained in $B_{i}$ either,
because by Remark 9.6(iv), if $m(R_{i}/(\xi_{0}))$ would be contained
in one of $m(R_{(0)}),\cdots,m(R_{(i-1)})$, then $R_{i}/(\xi_{0})$
would be a face of one of $R_{(0)},\cdots,R_{(i-1)}$, which is obviously
not the case. On the other hand, $\xi_{0}$ is the \emph{only} edge
of $R_{(i)}$ for which $m(R_{(i)}/(\xi_{0}))$ is not contained in
$B_{i}$ (indeed, the only other candidate would be $\xi_{i}$, but
$R_{(i)}/\xi_{i}=R_{(0)}$ and $m(R_{(0)})\subseteq B_{1}$). So,
we can form the pushout\[
\xymatrix{*++{m(\Lambda^{\xi_{0}}R_{(i)})}\ar@{>->}[d]\ar[r] & *++{B_{i}}\ar@{>->}[d]\\
m(R_{(i)})\ar[r] & B_{i+1}}
\]

Next, consider for each $i<j$ the tree \[
R_{(ij)}=R/(\xi_{1},\cdots,\hat{\xi}_{i},\cdots,\hat{\xi}_{j},\cdots,\xi_{n}),\]
and order these lexicographically, say as\[
R_{1}^{2},\cdots,R_{u}^{2},\,\,\,\,\,\,\,\,(u={{n \choose 2}}).\]
 We are going to form anodyne extensions of $B_{n+1}$ by using these
trees, \[
\xymatrix{B_{n+1}\,\,\ar@{>->}[r] & B_{n+2}\,\,\ar@{>->}[r] & \cdots\,\,\ar@{>->}[r] & B_{n+1+u}}
\textrm{, }\]
treating $R_{p}^{2}$ in the step to form $B_{n+p}\rightarrowtail B_{n+p+1}$
(for each $p=1,\cdots,u$). Suppose $B_{n+p}$ has been formed, and
consider $R_{p}^{2}=R_{(ij)}$ say. If $m(R_{p}^{2})\subseteq B_{n+p}$
then let $B_{n+p+1}=B_{n+p}$. If not, then surely $m(R_{p}^{2})\nsubseteq A_{k}\cup B$,
so by assumption (ii) $m(R_{p}^{2}/(\xi_{0}))=m(R/(\xi_{0},\xi_{1},\cdots,\hat{\xi}_{i},\cdots,\hat{\xi}_{j},\cdots,\xi_{n}))$
is not contained in $A_{k}\cup B$. On the other hand, Remark 9.6(iv)
implies that $m(R_{p}^{2}/(\xi_{0}))$ cannot be contained in any
of $m(R_{1}),\cdots,m(R_{n}),m(R_{1}^{2}),\cdots,m(R_{p-1}^{2})$
either. So $m(R_{p}^{2}/(\xi_{0}))$ is not contained in $B_{n+p}$.
As before, $\xi_{0}$ is the \emph{only} inner edge $\zeta$ for which
$m(R_{p}^{2}/(\zeta))$ is not contained in $B_{n+p}$. So we can
form the pushout\[
\xymatrix{*++{m(\Lambda^{\xi_{0}}(R_{(p)}^{2}))}\ar[r]\ar@{>->}[d] & *++{B_{n+p}}\ar@{>->}[d]\\
m(R_{(p)}^{2})\ar[r] & B_{n+p+1}}
\]

Next consider for each $i_{1}<i_{2}<i_{3}$ the tree \[
R_{(i_{1}i_{2}i_{3})}=R/(\xi_{1},\cdots,\hat{\xi}_{i_{1}},\cdots,\hat{\xi}_{i_{2}},\cdots,\hat{\xi}_{i_{3}},\cdots,\xi_{n})\]
 and adjoin the pushout along\[
m(\Lambda^{\xi_{0}}R_{(i_{1}i_{2}i_{3})})\rightarrowtail m(R_{(i_{1}i_{2}i_{3})})\]
 if necessary. Continuing in this way for all $l=0,1,\cdots,n-1$
and all sub-sequences $i_{1}<\cdots<i_{l}$ and corresponding trees
\[
R/(\xi_{1},\cdots,\hat{\xi}_{i_{1}},\cdots,\hat{\xi}_{i_{l}},\cdots,\xi_{n})\textrm{, }\]
 we end up with a sequence of anodyne extensions\[
B_{1}\rightarrowtail\cdots\rightarrowtail B_{q}\]
 where $q=\Sigma_{l=0}^{n-1}{{n \choose l}}=2^{n}-1$, and where $m(R/(\xi_{i}))$
is contained in $B_{q}$ for each $i=1,\cdots,n$. In the final step,
and exactly as before, we let $B_{2^{n}}=B_{2^{n}-1}$ if $m(R)\subseteq B_{2^{n}-1}$;
and if not, we form the pushout\[
\xymatrix{*++{m(\Lambda^{\xi_{0}}(R))}\ar[r]\ar@{>->}[d] & *++{B_{2^{n}-1}}\ar@{>->}[d]\\
m(R)\ar[r] & B_{2^{n}}}
\]
 Then $B_{2^{n}}$ is the pushout of $A_{0}\cup B$ and $m(R)$ over $(A_{0}\cup B)\cap m(R)$
(because every face $F$ of $R$ for which $m(F)$ is contained in $A_{0}\cup B$
occurs in some corner of the pushouts taken in the construction of
the $B_{i}$). This proves the lemma.
\end{proof}
Consider the tree $T_{k+1}$, and look at all lowest occurrences of
the $S$-colour $e$ (Recall $e$ is the fixed edge in $S$, occurring
in the statement of Proposition 9.2).  More precisely, let $e_{i}=(e,x_{i})$
for $i=1,\cdots,t$ be all the edges in $T_{k+1}$ whose $S$-colour
is $e$, while the $S$-colour of the edge immediately below it isn't.
This means that $(e,x_{i})$ is an edge having a white vertex at its
bottom. Let $\beta_{i}$ be the branch in $T_{k+1}$ from the root
to and including this edge $e_{i}$. Each such $\beta_{i}$ is an
initial segment in $T_{k+1}$, to which we will refer as the \emph{spine}
through $e_{i}$. For example, this is a picture of a spine in $T_{k+1}$,\[
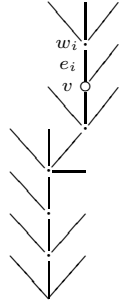
\begin{array}{ccc}
 & \quad & \xyC{10pt}\xyR{10pt}\xymatrix{ &  & *{}\ar@{-}[dr] & *{}\ar@{-}[d] & *{}\ar@{-}[dl]\\
 &  & *{}\ar@{}[r]|{w_{i}} & *{\cdot}\ar@{-}[d]_{e_{i}} & *{}\ar@{-}[dl]\\
 &  & *{}\ar@{-}[dr]\ar@{}[r]|{v} & *{\circ}\ar@{-}[d] & *{}\ar@{-}[dl]\\
 & *{}\ar@{-}[dr] & *{}\ar@{-}[d] & *{\cdot}\ar@{-}[dl]\\
 & *{}\ar@{-}[dr] & *{\cdot}\ar@{-}[d] & *{}\ar@{-}[l]\\
 & *{}\ar@{-}[dr] & *{\cdot}\ar@{-}[d] & *{}\ar@{-}[dl]\\
 & *{}\ar@{-}[dr] & *{\cdot}\ar@{-}[d] & *{}\ar@{-}[dl]\\
 &  & *{}\\
}
\end{array}\]
corresponding to the edge $e$ in $S$
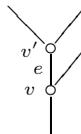
\[
\xymatrix{*{}\ar@{-}[dr] &  & *{}\ar@{-}[dl]\\
\ar@{}[r]|{ v'} & *{\circ}\ar@{-}[d]_{e} & *{}\ar@{-}[dl]\\
\ar@{}[r]|{ v} & *{\circ}\ar@{-}[d]\\
 & *{}}
\]

\begin{lem}
Let $R,Q_{1},\cdots,Q_{p}$ be initial segments in $T_{k+1}$, as
in the preceding lemma, and suppose condition (i) of that lemma is
satisfied. Then for any spine $\beta_{i}$ contained in $R$, the
edge $e_{i}\in\beta_{i}$ is characteristic for $R$ with respect
to $Q_{1},\cdots,Q_{p}$. 
\end{lem}
\begin{proof}
We have to check condition (ii) of Lemma 9.7. So, suppose $F$ is
an inner face of $R$, and suppose $m(F/(e_{i}))$ is contained in
$A_{k}\cup B=A_{0}\cup m(T_{1})\cup\cdots\cup m(T_{k})\cup m(Q_{1})\cup\cdots\cup m(Q_{p})\subseteq m(T_{k+1})$.
Since $m(F/(e_{i}))$ is isomorphic to the representable dendroidal
set $\Omega[F/(e_{i})]$, it must be contained in one of the dendroidal
sets constituting this union. But, if $m(F/(e_{i}))$ is contained
in $A_{0}$, then by Remark 9.6 (ii) $m(F)$ is also contained in
$A_{0}$. (the only colour occurring in $F$ but possibly not in $F/(e_{i})$
is the $S$-colour $e$). And, if $m(F/(e_{i}))$ is contained in
$m(T_{j})$ for some $j\le k$, then there must be a tensor product
relation applying to the image of $F/(e_{i})$, which allows a black
vertex to move up so as to get into an earlier $T_{j}$, as in:\[
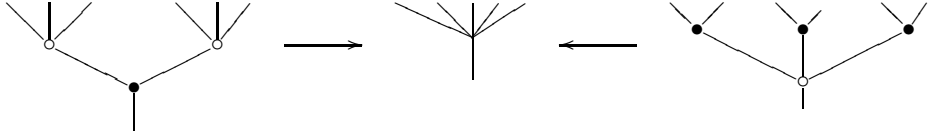
\begin{array}{ccccc}
\xyC{10pt}\xyR{10pt}\xymatrix{*{}\ar@{-}[dr] & *{}\ar@{-}[d] & *{}\ar@{-}[dl] &  & *{}\ar@{-}[dr] & *{}\ar@{-}[d] & *{}\ar@{-}[dl]\\
 & *{\circ}\ar@{-}[drr] &  &  &  & *{\circ}\ar@{-}[dll]\\
 &  &  & *{\bullet}\ar@{-}[d]\\
 &  &  & *{}}
 & \xymatrix{\\*{}\ar[rr] &  & *{}\\
}
 & \xymatrix{*{}\ar@{-}[drr] & *{}\ar@{-}[dr] & *{}\ar@{-}[d] & *{}\ar@{-}[dl] & *{}\ar@{-}[dll]\\
 &  & *{}\ar@{-}[d]\\
 &  & *{}}
 & \xymatrix{\\*{} &  & *{}\ar[ll]\\
}
 & \xyC{4pt}\xyR{4pt}\xymatrix{*{}\ar@{-}[dr] &  & *{}\ar@{-}[dl] &  & *{}\ar@{-}[dr] &  & {}\ar@{-}[dl] &  & *{}\ar@{-}[dr] &  & *{}\ar@{-}[dl]\\
 & *{\bullet}\ar@{-}[ddrrrr] &  &  &  & *{\bullet}\ar@{-}[dd] &  &  &  & *{\bullet}\ar@{-}[ddllll]\\
\\ &  &  &  &  & *{\circ}\ar@{-}[d]\\
 &  &  &  &  & *{}}
\end{array}\]
where the left tree is in $T_{k+1}$, the right tree is in $T_{j}$,
and the middle one in $F/(e_{i})$. 

But then the same relation must apply to the image of $F$, because
the edge $e_{i}$, having a white vertex at its root, cannot contribute
to this relation. Finally, if $m(F/(e_{i}))$ is contained in $m(Q_{l})$
for some $l\le p$, then by Remark 9.6 (iv), we have $m(R)\subseteq m(Q_{l})$.
So a fortiori, $m(F)$ is contained in $m(Q_{i})$. This proves the
lemma. 
\end{proof}
Recall that our aim is to prove for $A_{k}=A_{0}\cup m(T_{1})\cup\cdots\cup m(T_{k})$
that each inclusion\[
A_{k}\rightarrowtail A_{k+1}\]
 is anodyne. Consider the tree $T_{k+1}$, and let $\beta_{1},\cdots,\beta_{t}$
be all the spines contained in it. We shall prove by induction that
$A_{k}\rightarrowtail A_{k}\cup m(R_{1})\cup\cdots\cup m(R_{q})$
is anodyne, for any family $R_{1},\cdots,R_{q}$ of initial segments
each of which contains at least one such spine. The induction will
be on the number of such initial segments as well as on their size.
When applied to the maximal initial segment $T_{k+1}$ itself, this
will show that $A_{k}\rightarrowtail A_{k}\cup m(T_{k+1})=A_{k+1}$
is anodyne, as claimed. The precise form of induction is given by
the following lemma. 

\begin{lem}
Fix $l$ with $0\le l\le t$. Let $Q_{1},\cdots,Q_{p}$ be a family
of initial segments in $T_{k+1}$, each containing at least one and
at most $l$ spines. Let $R_{1},\cdots,R_{q}$ be initial segments,
each of which contains $l+1$ spines. Then the inclusion \[
A_{k}\rightarrowtail A_{k}\cup B\cup C\]
 for $B=m(Q_{1})\cup\cdots\cup m(Q_{p})$ and $C=m(R_{1})\cup\cdots\cup m(R_{q})$,
is anodyne.
\end{lem}
\begin{proof}
We can measure the size of each of the initial segments $R_{j}$ by
counting the number of vertices in $R_{j}$ which are not on one of
the $l+1$ spines. If this number is not bigger than $u$, we say
that $R_{j}$ has size at most $u$, and write $size(R_{j})\le u$.
Let $\Lambda(l,u)$ be the assertion that the lemma holds for $l$,
for any families $\{ Q_{i}\}$ and $\{ R_{j}\}$ where the $R_{j}$
all have $size(R_{j})\le u$. We will prove $\Lambda(l,u)$ by induction,
first on $l$ and then on $u$.\emph{}\\

\emph{Case $l=0$:} This is the case where there are no $Q$'s, i.e.,
$p=0$. For $l=0$, first consider the case where $u=0$ also. Then
each $R_{i}$ is itself a spine, say $\beta_{i}$, with top inner
edge $e_{i}$ running from a copy of $v$ to a copy $w_{i}$ of $w$.
We will prove that each of the inclusions\[
A_{k}\cup m(R_{1})\cup\cdots\cup m(R_{i-1})\rightarrowtail A_{k}\cup m(R_{1})\cup\cdots\cup m(R_{i})\]
 for $i=0,\cdots,q$, is anodyne. If $R_{i}=\beta_{i}$ coincides
with one of the earlier spines $R_{j},\, j<i$, then there is nothing
to prove. If $R_{i}$ is a different spine, then its outer top face
is contained in $A_{0}$ because it misses the vertex $v'$ which
is above $e$ in $S$. So condition (i) of Lemma 9.7 is satisfied,
where $R_{i},R_{1},\cdots,R_{i-1}$ take the role of $R,Q_{1},\cdots,Q_{p}$
in that lemma. By Lemma 9.8, the edge $e_{i}\in R_{i}$ is characteristic
so Lemma 9.7 gives that $A_{k}\cup m(R_{1})\cup\cdots\cup m(R_{i-1})\rightarrowtail A_{k}\cup m(R_{1})\cup\cdots\cup m(R_{i})$
is anodyne, as claimed. The composition of these inclusions will then
be anodyne also, which proves the statement $\Lambda(0,0)$.\\

Suppose now that $\Lambda(0,u)$ has been proved, and consider families
$R_{1},\cdots,R_{q}$ of initial segments which are each of size not
bigger than $u+1$. Suppose that among these, $R_{1},\cdots,R_{q'}$
actually have size not bigger than $u$, while $R_{q'+1},\cdots R_{q}$
have size $u+1$. We shall prove that \[
A_{k}\rightarrowtail A_{k}\cup m(R_{1})\cup\cdots\cup m(R_{q})\]
 is anodyne, by induction on the number $r=q-q'$ of initial segments
that have size $u+1$. If $r=0$, this holds by $\Lambda(0,u)$. Suppose
we have proved this for \emph{any} family with not more than $r$
initial segments of size $u+1$, and consider such a family $R_{1},\cdots,R_{q}$
where $q-q'=r+1$. Write $\beta_{q}$ for the spine contained in $R_{q}$
(there is only one such because we are still in the case $l=0$).
For a top outer face $\partial_{x}(R_{q})$ of $R_{q}$, \emph{either}
$x\ne w_{q}$ so that $\partial_{x}(R_{q})$ still contains $\beta_{q}$
but has size at most $u$, \emph{or} $x=w_{q}$ so that $m(\partial_{x}(R_{i}))$
is contained in $A_{0}$ because it misses the vertex $v'$ immediately
above $e$ in $S$. Thus, if we let \[
P=m(R_{1})\cup\cdots\cup m(R_{q-1})\cup\bigcup_{x}m(\partial_{x}(R_{q}))\]
 where $x$ ranges over all the top vertices in $R_{q}$, then by
the fact that $\Lambda(0,u+1)$ is assumed to hold for $r=(q-1)-q'$,
\[
\begin{array}{ccc}
\quad\quad\quad\quad\quad\quad\quad\quad\quad\quad\quad\quad\quad & A_{k}\rightarrowtail A_{k}\cup P & \quad\quad\quad\quad\quad\quad\quad\quad\quad\quad\quad\quad\quad(1)\end{array}\]
is anodyne. To prove that $A_{k}\cup P\rightarrowtail A_{k}\cup P\cup m(R_{q})=A_{k}\cup m(R_{1})\cup\cdots\cup m(R_{q})$
is anodyne as well, we can now apply Lemma 9.7. Indeed, the family
of initial segments containing $P$ is made to contain the images
of all the top faces of $R_{q}$, and $e_{q}\in R_{q}$ is characteristic
by Lemma 9.8. This proves that $A_{k}\cup P\rightarrowtail A_{k}\cup P\cup m(R_{q})$
is anodyne, as claimed. When composed with (1), we find that $A_{k}\rightarrowtail A_{k}\cup m(R_{1})\cup\cdots\cup m(R_{q})$
is anodyne. This proves $\Lambda(0,u+1)$ and completes the inductive
proof of $\Lambda(l,u)$ for $l=0$ and all $u$. \\

Suppose now that we have proved $\Lambda(l',u)$ for all $l'\le l$
and all $u$. We will now prove $\Lambda(l+1,u)$ by induction on
$u$. 

\emph{Case} $u=0$: This is the assertion that for any given initial
segments \[
Q_{1},\cdots,Q_{p},R_{1},\cdots,R_{q}\]
 of $T_{k+1}$, where the $Q_{j}$ contain at most $l$ spines while
each $R_{i}$ is made up out of exactly $l+1$ spines (and no other
vertices), the inclusion \[
\begin{array}{ccc}
\quad\quad\quad & A_{k}\rightarrowtail A_{k}\cup m(Q_{1})\cup\cdots\cup m(Q_{p})\cup m(R_{1})\cup\cdots\cup m(R_{q}) & \quad\quad\quad(2)\end{array}\]
is anodyne. We shall prove by induction on $q$ that this holds for
all $p$. For $q=0$, the conclusion follows by the inductive assumption
that $\Lambda(l,u)$ holds. Suppose the assertion holds for $q-1$,
and consider $R_{q}$. Each top vertex of $R_{q}$ lies at the end
of a spine, so $\partial^{x}(R_{q})$ contains at most $l$ spines.
Let \[
D=m(Q_{1})\cup\cdots\cup m(Q_{p})\cup\bigcup_{x}\partial_{x}(R_{q})\]
where $x$ ranges over the top vertices of $R_{q}$. Then, by the
assumption for $q-1$, \[
\begin{array}{ccc}
\quad\quad\quad\quad\quad\quad\quad & A_{k}\rightarrowtail A_{k}\cup D\cup m(R_{1})\cup\cdots\cup m(R_{q-1}) & \quad\quad\quad\quad\quad\quad\quad(3)\end{array}\]
is anodyne. To prove that (2) is anodyne, it then suffices to apply
Lemma 9.7, and show that $R_{q}$ has a characteristic edge with respect
to the family of initial segments containing the union $D\cup m(R_{1})\cup\cdots\cup m(R_{q-1})$
in (3). But by Lemma 9.8, any top edge $e_{q}$ of $R_{q}$ is characteristic.
This proves $\Lambda(l+1,u)$, for $u=0$.

\emph{Case} $u+1$: Suppose now $\Lambda(l+1,u)$ holds. To prove
$\Lambda(l+1,u+1)$, consider families\[
\begin{array}{ccc}
\quad\quad\quad\quad\quad\quad\quad & Q_{1},\cdots,Q_{p},R_{1},\cdots,R_{q'},R_{q'+1},\cdots,R_{q} & \quad\quad\quad\quad\quad\quad\quad(4)\end{array}\]
of initial segments in $T_{k+1}$, where the $Q_{i}$ contain at most
$l$ spines, the $R_{i}$ contain exactly $l+1$ spines, the $R_{1},\cdots,R_{q'}$
are of size not more than $u$, and $R_{q'+1},\cdots,R_{q}$ are of
size exactly $u+1$. We will show by induction on the last number
$r=q-q'$ that for any such family, the inclusion\[
\begin{array}{ccc}
\quad\quad\quad & A_{k}\rightarrowtail A_{k}\cup m(Q_{1})\cup\cdots\cup m(Q_{p})\cup m(R_{1})\cup\cdots\cup m(R_{q}) & \quad\quad\quad(5)\end{array}\]
is anodyne. For $r=q-q'=0$ there is nothing to prove, because this
is the case covered by $\Lambda(l+1,u)$. Suppose we have proved that
(5) is anodyne for \emph{any} family (4) with $q-q'\le r$, and consider
such a family with $q-q'=r+1$. The initial segment $R_{q}$ has two
kinds of top outer faces, namely the $\partial_{x}(R_{q})$ which
remove the top of a spine, and the $\partial_{x}(R_{q})$ where $x$
does not lie on a spine. Outer faces of the first kind contain $l$
spines only, and outer faces of the second kind are of size not more
than $u$. Let \[
D=m(Q_{1})\cup\cdots\cup m(Q_{p})\cup\bigcup_{x}m(\partial_{x}R_{q})\]
 where $x$ ranges over the top vertices of $R_{q}$ which are on a
spine. Let\[
E=m(R_{1})\cup\cdots\cup m(R_{q'})\cup\bigcup_{x}m(\partial_{x}(R_{q}))\]
 where $x$ ranges over the top vertices of $R_{q}$ which are not on a
spine. Then, by the assumption that $\Lambda(l+1,u+1)$ has been established
for families (4) where $q-q'\le r$, we see that \[
\begin{array}{ccc}
\quad\quad\quad\quad\quad\quad\quad & A_{k}\rightarrowtail A_{k}\cup D\cup E\cup R_{q'+1}\cup\cdots\cup R_{q-1} & \quad\quad\quad\quad\quad\quad\quad(6)\end{array}\]
is anodyne. The union $D\cup E$ is made to contain all the images
$\partial_{x}(R_{q})$ of top faces $\partial_{x}(R_{q})$ of $R_{q}$,
and by Lemma 9.8, any edge $e_{q}$ on the top of a spine $\beta_{q}$
in $R_{q}$ is characteristic with respect to the family of initial
segments making up the union on the right-hand-side of (6). So by
Lemma 9.7, the map \[
A_{k}\cup D\cup E\cup R_{q'+1}\cup\cdots\cup R_{q-1}\rightarrowtail A_{k}\cup D\cup E\cup R_{q'+1}\cup\cdots\cup R_{q}\]
 is anodyne. When composed with (6), this gives (5), and proves the
case $u+1$. 

This established $\Lambda(l+1,u+1)$ and completes, for $l+1$, the
induction on $u$, thus completing the proof. 
\end{proof}

\end{document}